\documentclass[a4paper,english]{article}
\title{The CS decomposition and conditional negative correlation inequalities for  determinantal processes.} 
\author{Andr\'e Goldman}

\usepackage[english]{babel}  
\usepackage{indentfirst}  
                                               
\usepackage[latin1]{inputenc}  
\usepackage[T1]{fontenc}  
\usepackage{amssymb,amsmath,amscd}

\newtheorem{notation}{Notation}
\newtheorem{lemma}{Lemma}

\newtheorem{prop}[lemma]{Proposition}
\newtheorem{theorem}{Theorem}
\newtheorem{remark}{Remark}
\newtheorem{example}{Example}

\newenvironment{proof}[1][]%
{\noindent{\bf Proof{}#1:}\par\nobreak}{\nopagebreak\hspace*{\fill} $\square$\vspace*{1.4ex}\par}

\begin{document}

\maketitle

\begin{abstract}For a conditional process of the form
$(\phi \vert A_{i} \not\subset  \phi)$ where $\phi$ is a 
determinantal process we obtain a new negative correlation inequalities. Our approach relies upon the underlying geometric structure of the elementary discrete determinantal  
processes by using the canonical representation of a pair of subspaces in terms of principal vectors and angles, as well as the classical CS decomposition.

\end{abstract}
\section{Introduction and statement of results.}\label{S1}
Determinantal point processes (DPP) are probabilistic models that arise in a wide range of theoretical and applied areas.
Let us cite,  to mention a few recent applied works, 
\cite{ref2},  \cite{ref9}, \cite{ref10}, \cite{ref18}, \cite{ref23}, \cite{ref24}.   
From the theoretical point of view,  determinantal point processes could be defined in the general  locally compact Polish spaces setting, as point processes associated to some locally square integrable, Hermitian, positive semidefinite, locally trace-class operators; a good overview of the main conceptual basis and properties can be found  in \cite{ref20} and in the bibliography therein; see also \cite{ref16}. \\
There exists also various extensions and variants as for example L-ensemble \cite{Bor}, extended L-ensembles \cite{Tr}, strongly Rayleigh point processes \cite{ref27} and so on. A common feature of all these  processes 
is to  exhibit repulsion between points and to offer, among others,  efficient algorithms for sampling and conditioning. Despite of their widespread use the study of their properties is, in general, not easy and several questions and conjectures are still unsolved (even in the most elementary setting). Notice also that some results for the most general processes can be infered  \cite{ref11}, \cite{ref20}, \cite{ref21},  
from the corresponding results of the basic processes.\\
Recall that the basic elementary determinantal point process can be described via the exterior product concept, as follows.\\
Fix $1< p < N$ and let 
  $ \mathfrak{A}=\{a^{1},...,a^{p}\} $, $ 1<p<N $, be a set of orthonormal vectors
  in $ \mathbb{R}^{N} $. Denote 
 
$$
a^{i}=(a^{i}_{1},\dots,a^{i}_{N})^{t}, i=1,\dots,p
$$
and
$$
a_{i}=(a^{1}_{i},\dots,a^{p}_{i}), i=1,\dots,N.
$$
The associated determinantal process $ \phi(\mathfrak{A}) $
is a point process, view as a random subset of $\mathcal{N}=\{1,\dots,N\}$
of cardinality $ \vert \phi(\mathfrak{A})\vert = p $, characterized \cite{ref19}, \cite{ref20}
by the formulas
\begin{equation}\label{I-1} 
\mathbb{P}\{\{i_{1},\dots,i_{p}\} = \phi\}=
 \mid(\bigwedge_{i=1}^{p}a^{i})_{\{i_{1},\dots,i_{p}\}}\mid^{2}
= [\det((a_{i_{j}}^{k})_{k,j=1,\dots,p})]^{2}
\end{equation}
for all subsets $\{i_{1},\dots,i_{p}\} \subset \mathcal{N}$.
Note also that (\ref{I-1}) implies 
\begin{equation}\label{I-2}
\mathbb{P}\{\{i_{1},\dots,i_{k}\} \subset \phi\}=
\parallel \bigwedge_{j=1}^{k}a_{i_{j}}\parallel ^{2}
\end{equation}
for all $1\leq k \leq p$.

Let $ E=E(\mathfrak{A})\subset
 \mathbb{R}^{N}$ 
 be 
the vector space spanned by  $ \mathfrak{A} $. For all sets of linearly independent vectors $v^{i}\in E$, $i=1,\dots,p$,   we have $\bigwedge_{i=1}^{p}v^{i}=a\bigwedge_{i=1}^{p}a^{i}$ with $a\neq 0$ thus, in particular, if 
$ \tilde{\mathfrak{A}}=\{\tilde{a}^{1},...,\tilde{a}^{p}\} $
is another orthonormal basis of $ E=E(\mathfrak{A}) $ then  
$$
\mid(\bigwedge_{i=1}^{p}a^{i})_{\{i_{1},\dots,i_{p}\}}\mid
=
\mid(\bigwedge_{i=1}^{p}\tilde{a}^{i})_{\{i_{1},\dots,i_{p}\}}\mid
$$ 
for every
$\{i_{1},\dots,i_{p}\}\subset \mathcal{N}$ and consequently
$ \phi(\mathfrak{A})=\phi(\tilde{\mathfrak{A}}) $. \\
Remark also that if $ \mathfrak{A}^{\perp}=\{a^{p+1},...,a^{N}\} $ is
an orthonormal basis of  the orthogonal complement 
$ E(\mathfrak{A})^{\perp} $ of  $E(\mathfrak{A})$ in $\mathbb{R}^{N}$
then obviously
$$
\phi(\mathfrak{A}^{\perp})=\mathcal{N}\setminus \phi(\mathfrak{A}).
$$
An example of a non-trivial basic determinantal process is given by uniform spanning tree measure on a finite connected graph G. Roughly speaking, if G fixed, is arbitrary edge-oriented and M is the vertex-edge incidence matrix (the columns been indexed by vertexes) then the determinantal process associated to the vector space spanned by column vectors but one, provides a uniform probability on spanning trees. This result due \cite{ref7} is called the Transfer Current Theorem. For more details see \cite{ref20}. 
Some extensions of this result are given in \cite{ref4} with a serie of open questions and conjectures.\\
The repulsion between points of the process $ \phi(\mathfrak{A}) $  is reflected by its negative dependence properties \cite{ref2}, \cite{KuTa}, \cite{ref23}, \cite{ref26}, . The simplest one is  the negative correlation inequality 
\begin{equation}\label{I-3}
P \left\lbrace A \subset  \phi \vert B\subset  \phi
\right\rbrace\leq P \left\lbrace A \subset  \phi\right\rbrace  
\end{equation} 
where $A$ and $B$ are disjoint subsets  of $\mathcal{N}$ such that 
$\mathbb{P}\lbrace B  \subset \phi \rbrace > 0$.
Note that by formula (\ref{I-2}) the inequality (\ref{I-3}) can be rephrased as 
\begin{equation}\label{I-4}
\parallel \bigwedge_{i \in A\cup B}a_{i}\parallel ^{2}
\leq \parallel \bigwedge_{i \in A}a_{i}\parallel ^{2} \times
\parallel \bigwedge_{i \in B}a_{i}\parallel ^{2}
\end{equation}
which is the well-known Hadamard-Whithney inequality.\\
\begin{remark}\label{R3}
For a finite set $E\subset \mathbb{R}^{p}$ of vectors the formula $\bigwedge_{x\in E}x$ designates exterior product of vectors
$x\in E$ when these are arbitrarily numbered. For  different 
numberings the associated  exterior products may differ (only) by sign which in our context doesn't matter.
\end{remark}
 A  strengthening  of (\ref{I-3}) is the (less obvious) association inequality
 obtained by R.Lyons in \cite{ref19}
\begin{equation}\label{I-5}
P\left\lbrace \phi \in \mathcal{A} \vert \phi \in \mathcal{B}\right\rbrace 
\leq P\left\lbrace \phi \in \mathcal{A} \right\rbrace  
\end{equation}
where $\mathcal{A} $ and $\mathcal{B}$ are increasing events which are measurable wih complementary subsets and $P(\phi \in \mathcal{B} )> 0 $. An event $\mathcal{A} \subset 2^{\mathcal{N}}$, 
 is called increasing if whenever $A\in \mathcal{A} $ and  $n \in \mathcal{N}$, we have also $A\cup \{n\}\in \mathcal{A} $.\\
 The next step is the BK inequality. For a pair $\mathcal{A} $, $\mathcal{B} $ $\subset 2^{\mathcal{N}}$ of increasing events the disjoint intersection  $\mathcal{A} \circ \mathcal{B}$ is  defined \cite{ref6} by
\begin{equation}\label{I-6}
\mathcal{A}\circ\mathcal{B} 
=\{ K \subset \mathcal{N}:\exists \quad L\in \mathcal{A}, \quad M \in \mathcal{B}, \quad L,M \neq \emptyset, \quad L \cap M = \emptyset,
K \supset L\cup M
 \}.
 \end{equation}
A point process $\psi$ on $\mathcal{N}$ is said to have the van den Berg - Kesten property (in short the BK property) if 
\begin{equation}\label{I-7}
\mathbb{P}\{ \psi \in \mathcal{A}\circ\mathcal{B} \vert \psi \in  \mathcal{B}\}
\leq 
\mathbb{P}\{\psi \in  \mathcal{A}\}
\end{equation}

for every pair of increasing events such that $\mathbb{P}\{\psi \in  \mathcal{B}\}>0$. In \cite{ref6} (see also \cite{ref5}) J. van den Berg and  H.Kesten proved that inequality (\ref{I-7}) is satisfied when $\psi$ is related to a product probability on $2^{\mathcal{N}}$. In the basic determinantal process setting, the Conjecture 4.6. of \cite{ref4} which  states that the same is true for the spanning trees determinantal point processes is still unsolved. The question of whether general determinantal processes have the BK property was raised in \cite{ref19}. Note also
that for increasing events $\mathcal{A}$ and $\mathcal{B}$   which are measurable with complementary subsets we have 
$\mathcal{A} \circ \mathcal{B} = \mathcal{A} \cap \mathcal{B}$ and thus  (\ref{I-7}) implies (\ref{I-5}) .\\
In \cite{ref12} it it was conjectured that for all basic determinantal processes the following conditional inequalities
hold:  
 \begin{equation}\label{I-8}
\begin{split}
& \mathbb{P}\lbrace  A \not \subset \phi \mid A_{i} \not \subset \phi, \forall i=1,\dots,n\rbrace \\
& \leq \mathbb{P}\lbrace A \not \subset \phi \mid A_{i} \not \subset \phi, \forall i=1,\dots,n-1\rbrace
 \end{split} 
\end{equation}
for all disjoint subsets $A$, $A_{i}$, $i=1,\dots,n$,   of $\mathcal{N}$
 such that  $\mathbb{P}\lbrace  A_{i} \not \subset \phi, \forall i=1,\dots , n\rbrace  > 0 $. 
 The inequalities (\ref{I-8}) are closely related to the BK inequality namely
they imply (theorem 2 in \cite{ref12}) that the BK inequality is satisfied for increasing events $\mathcal{A}$ and $\mathcal{B}$ which are generated by disjoint sets. When all the sets above are  reduced to being  simple points then  the  inequalities (\ref{I-8}) are satisfied. This follows easily from the fact that in this particular case the conditioned processes $(\phi \mid x_{i}\notin \phi \quad i=1,\dots,n)$ are still determinantal (see for example point 4 of proposition 1 in \cite{ref12}). However, this later doesn't hold for general sets as the following example shows.
\begin{example}\label{ex} 
Let 
  $ \mathfrak{A}=\{a^{1},a^{2}\} $  be orthonormal vectors 
  in $ \mathbb{R}^{4} $  
 and   $ \phi = \phi(\mathfrak{A}) $ the associated  determinantal process. We suppose that
 \begin{enumerate}
 \item
 $0 < P\{\phi \neq \{1,2\}\}=\kappa < 1$
 \item
$P\{\phi = \{1,4\}\}>0$
\item
$P\{\phi = \{2,3\}\}>0$,
\end{enumerate}
then the proces  $ \psi = (\phi \vert\{1,2\} \not \subset \phi) $ 
, $\vert \psi \vert = 2$
satisfies
\begin{enumerate}
\item[4] 
 $P\{\psi = \{1,i\}\} = P\{\phi = \{1,i\}\}/\kappa $
 $i=3,4$  
 \item[5]
 $P\{\psi = \{2,i\}\} = P\{\phi = \{2,i\}\}/\kappa $
 $i=3,4$. 
\end{enumerate}
It follows from  2-3 that $P\{ i \in \psi\} > 0$, $i=1,2,3,4$.
Now suppose that $\psi$ is determinantal. Then there exist 
orthonormal vectors  $ \mathfrak{\tilde A}=\{\tilde a^{1},\tilde a^{2}\} $ in $ \mathbb{R}^{4} $ such that
$\psi = \psi (\mathfrak{\tilde A})$.
\\
But $\parallel \tilde a_{1}\wedge \tilde a_{2} \parallel^{2} = P\{\psi = \{1,2\}\}=0$ and
$\parallel \tilde a_{i} \parallel^{2}= P\{ i \in \psi\} > 0$, $i=1,2,$
 therefore there
exist $\alpha \neq 0 $ such 
that $\tilde a_{1}=\alpha \tilde a_{2}$ which implies that
$$P\{\psi = \{1,i\}\} = \parallel \tilde a_{1}\wedge \tilde a_{i} \parallel^{2} =
\alpha^{2} \parallel \tilde a_{2}\wedge \tilde a_{i} \parallel^{2}= \alpha^{2} P\{\psi = \{2,i\}\}, i=3,4 $$
and  from 4-5
\begin{equation}\label{I-9}
P\{\phi = \{1,i\}\} = \alpha^{2} P\{\phi = \{2,i\}\}, i=3,4. 
\end{equation}
In general this is not possible.
It suffices to take 
 $a^{1}=(1/\sqrt{2},0,1/\sqrt{2},0)^{t}$ and
$a^{2}=(0,1/\sqrt{2},0,1/\sqrt{2})^{t}$ to see that  the associated process $\phi_{0}$
does not fulfil (\ref{I-9}). 
\end{example}
\begin{remark}\label{R1}
A generalization
 of determinantal point processes is given by  strongly Rayleigh point processes which are introduced in the work of \cite{ref27}. A point process $\psi$ is said strongly Rayleigh if its generating polynomial 
$$g_{\psi}(z_{1},\dots,z_{N})=\sum_{S\subset \mathcal{N}}
\mathbb{P}\{ \psi=S\}\prod_{i\in S}z_{i}, \quad z_{i}\in \mathbb{C},\quad  i=1,\dots,N,$$
is stable, that is $g_{\psi}(z_{1},\dots,z_{N})\neq 0 $
whenever $\mathfrak{Im}(z_{i})>0$ 
for all $1 \leq i  \leq N$. It was proved in \cite{ref27}   that   strongly Rayleigh processes
has negative association and  that determinantal processes are strongly Rayleigh. 
As regards to example \ref{ex} note that the generating polynomial 
$g_{\psi_{0}}(z_{1},\dots,z_{4})=\frac{1}{3}(z_{1}z_{4}+z_{2}z_{3}+z_{3}z_{4})$
of the the conditioned process  $ \psi_{0} = (\phi_{0} \vert \{1,2\} \not \subset \phi_{0}) $ is not stable and consequently  $ \psi_{0}$ is not strongly Rayleigh.
 \end{remark}
 The purpose of this work is to shown that the inequalities (\ref{I-8})  are satisfied for n=2 and n=3. To this aim we reformulate  (\ref{I-8}) (which is an elementary exercice) in the equivalent form:
\begin{equation}\label{I-10}
\begin{split}
&\mathbb{P}\lbrace  A  \subset \phi \mid A_{i} \not \subset \phi, \forall i=1,\dots,n\rbrace \times \mathbb{P}\lbrace  B  \subset \phi \mid A_{i} \not \subset \phi, \forall i=1,\dots,n\rbrace \\
& -  \mathbb{P}\lbrace  A  \subset \phi, B  \subset \phi \mid A_{i} \not \subset \phi, \forall i=1,\dots,n\rbrace \geq 0.
 \end{split} 
\end{equation}
 Next, we will shown in Section \ref{s2} that in order to obtain (\ref{I-10}) it suffices to consider the case when A and B are simple points. In  this case we are able to compute  for $n=1$ and $n=2$ (which corresponds to $n=2$ and $n=3$ for  (\ref{I-8}))  the exact value  of the difference appearing in  (\ref{I-10}) in terms of some geometric characteristics of the vector space $E(\mathfrak{A})$. The proof relies upon the classical canonical representation of a pair of subspaces in terms of principal vectors and angles \cite{ref1}, \cite{ref13},
 \cite{ref17} and the  link etablished in \cite{ref12} between the classical CS decomposition (CSD) \cite{ref25} of partitioned unitary matrix and the basic determinantal processes which allows  to obtain a pertinent bases of the vector spaces $ E=E(\mathfrak{A})\subset
 \mathbb{R}^{N}$ and $ E^{\perp} =(\mathfrak{A})^{\perp} $. To clarify the latter: if we choose  $N-p - n +2>0$, $p=\vert \phi(\mathfrak{A})\vert$, and reorder the coordinates of the vectors in $\mathbb{R}^{N}$ with the order $(2,\dots,n-1,1,n,n+1,\dots,N)$ then the  first p columns (resp. the last N-p columns)
of the matrix below is an orthonormal basis of $E$ (resp. $ E^{\perp}$).
\begin{equation}\label{I-11}
\begin{split}
\begin{matrix}
2\\
\vdots\\
n-1\\
1\\
n\\
\vdots\\
N
\end{matrix} 
& \left[ \begin{matrix}
\cos \theta_{1}u^{1}_{1}  \dots  \cos \theta_{n-2}u^{n-2}_{1}
  &  0 \quad \dots \quad 0 \\
\vdots  \qquad & \vdots    \\
 \cos \theta_{1}u^{1}_{n-2}  \dots  \cos \theta_{n-2}u^{n-2}_{n-2}
  &   0 \quad \dots \quad 0 \\
\sin\theta_{1}V^{1}_{1}   \dots  \sin \theta_{n-2}V^{n-2}_{1} & W^{1}_{1} \dots W^{p-n+2}_{1}\\
\sin\theta_{1}V^{1}_{2}   \dots  \sin \theta_{n-2}V^{n-2}_{2} & W^{1}_{2} \dots W^{p-n+2}_{2}\\
\vdots  \qquad & \vdots    \\
\sin\theta_{1}V^{1}_{N-n}   \dots  \sin \theta_{n-2}V^{n-2}_{N-n}  & W^{1}_{N-n} \dots W^{p-n+2}_{N-n}
\end{matrix} \right. 
\\
& \\
&
\qquad   
 \left. \begin{matrix}
\sin \theta_{1}u^{1}_{1}  \dots  \sin \theta_{n-2}u^{n-2}_{1}
  &  0 \quad \dots \quad 0 \\
\vdots  \qquad & \vdots    \\
 \sin \theta_{1}u^{1}_{n-2}  \dots  \sin \theta_{n-2}u^{n-2}_{n-2}
  &  0 \quad \dots \quad 0 \\
-\cos\theta_{1}V^{1}_{1}   \dots  -\cos \theta_{n-2}V^{n-2}_{1} & \tilde{W}^{1}_{1}  \dots \  \tilde{W}^{N+2-n-p}_{1}\\
-\cos\theta_{1}V^{1}_{2}   \dots  -\cos \theta_{n-2}V^{n-2}_{2} & \tilde{W}^{1}_{2}  \dots   \tilde{W}^{N+2-n-p}_{2}\\
\vdots  \qquad & \vdots    \\
-\cos\theta_{1}V^{1}_{N-n}   \dots  -\cos \theta_{n-2}V^{n-2}_{N-n} & \tilde{W}^{1}_{N-n}  \dots   \tilde{W}^{N+2-n-p}_{N-n}
\end{matrix} \right]   
\end{split} 
\end{equation}    
Here, $ u^{1},\dots,u^{n} $ are orthonormal vectors in
$\mathbb{R}^{n-2}$; 
$ \mathfrak{V}=\{V^{1},\dots,V^{n}\} $,
$\mathfrak{W}= \{W^{1},\dots,W^{p-n}\} $, 
$ \mathfrak{\tilde{W}}=\{\tilde{W}^{1},\dots,\tilde{W}^{N-p-n}\} $,
are mutually orthonormal vectors in $\mathbb{R}^{N-n+2}$;
the angles appearing in the matrix above are principal Jordan angles between
the space E and the basic subspace
$$
\mathbb{R}_{\{2,\dots,n-2\}}^{N}=\{x=(x_{k})\in\mathbb{R}^{N};x_{k}=0 \quad
if \quad k\notin \{2,\dots,n-2\}\};
$$
for more details see Section 2 of \cite{ref12}.\\
By making use of these particular bases it is shown (Proposition 2 in  \cite{ref12})
 that 
\begin{equation}\label{I-12}
\begin{split}
&  P\{1 \in \psi \}\times P\{ n \in \psi \}
- P\{1\in \psi, n \in \psi \}
\\
& = \lambda\times (\langle (V_{1},W_{1}),(V_{2},W_{2})  \rangle)^{2}\\
& = \lambda\times
\left(\langle a_{1},a_{n} \ \rangle +
\sum_{k=1}^{n-1}(-1)^{k}\sum_{2\leq i_{1}<\dots < i_{k}\leq n-1}
  \langle a_{1}\wedge(\bigwedge_{j=1}^{k} a_{i_{j})},
 a_{n}\wedge(\bigwedge_{j=1}^{k} a_{i_{j}}) \ \rangle
\right)^{2} 
\end{split}
\end{equation}
where $\psi = ( \phi\vert i \not\in \phi, i=2,\dots ,n-1)$,
and
$P\{i\not\in \phi, i=2,\dots ,n-1\}
=\prod_{i=1}^{n-2}\sin^{2}\theta_{i}  = \lambda^{-\frac{1}{2}} >0$ . 
\begin{remark}\label{R2}
Note that the values of the right and left sides of (\ref{I-12}) does not depend of the choice of the basis $ \mathfrak{A}=\{a^{1},...,a^{p}\} $ and that by reordering the indices  the  procedure described above works for every subset  $\{i_{2},\dots,i_{n-1}\}\subset  \mathcal{N}$.
\end{remark}
Now, for disjoint  subsets $\{i\}$, $\{j\}$, $A=\{i_{1},\dots,i_{n-2}\}$ of $\mathcal{N}$, 
$0<P\{A\in \phi\} < 1$, $\vert A \vert= n-2 \leq p$ we choose a basis as above (the  n-2 terms of  A taken at first and $i,j$ follow). In addition we note
\begin{equation}\label{I-13}
\tilde{v}_{1}=(V_{1}^{1}\sin\theta_{1},\dots V_{1}^{n-2}\sin\theta_{k})
\end{equation}
and
\begin{equation}\label{I-14}
\tilde{v}_{2}=(V_{2}^{1}\sin\theta_{1},\dots V_{2}^{n-2}\sin\theta_{k})
\end{equation}
such that
\begin{equation}\label{I-15}
a_{i}=(\tilde{v}_{1},W_{1}), \quad a_{j}=(\tilde{v}_{2},W_{2}).
\end{equation}
With these notation we prove  in Section \ref{s3} that :
\begin{theorem}\label{th1}
\item[1.]
If   $\vert A\vert =n-2=p$ then
\begin{equation}\label{I-16}
\begin{split} 
& P\{ i\in \phi \vert A \not \subset \phi\}
\times 
P\{ j \in \phi \vert A \not \subset \phi \} - P\{ i\in \phi, j \in \phi \vert A \not \subset \phi \} 
\\
& 
= (1-\kappa)^{-2} \times \big( \langle a_{i}, a_{j}\rangle^{2}  +  
\kappa\times (a_{i}\wedge a_{j})^{2} \big).
\end{split}
\end{equation}

\item[2.]
If $\vert A\vert =n-2<p$ then
\begin{equation}\label{I-17}
\begin{split}
& P\{ i\in \phi \vert A \not \subset \phi\}
\times 
P\{ j \in \phi \vert A \not \subset \phi \} - P\{ i\in \phi, j \in \phi \vert A \not \subset \phi \} 
\\
& = (1-  \kappa)^{-2}
\times \bigg[ \bigg(\langle \tilde{v}_{1},\tilde{v}_{2}\rangle +
(1-  \kappa)\langle W_{1},W_{2}\rangle \bigg)^{2} 
+  \kappa\parallel\tilde{v}_{1}\wedge \tilde{v}_{2}\parallel^{2}
\bigg]
\end{split}
\end{equation}
where $\kappa=P\{ A\subset \phi\}=(\prod_{i=1}^{p}\cos\theta_{i})^{2}$.
Moreover  we have
\begin{equation}\label{I-18}
\langle \tilde{v}_{1},\tilde{v}_{2}\rangle +
(1-  \kappa)\langle W_{1},W_{2}\rangle = \langle a_{i}, a_{j}\rangle
- \langle a_{i}\wedge(\bigwedge_{k=1}^{n-2} a_{i_{k}}), a_{j}\wedge(\bigwedge_{k=1}^{n-2} a_{i_{k})}\rangle.
\end{equation}
Note that the right sides of (\ref{I-16}) and (\ref{I-18}) do not depend of the choice of the basis.\\
 As a consequence, the inequality
\begin{equation}
  \mathbb{P}\lbrace  A \not \subset \phi \mid A_{i} \not \subset \phi,   i=1,2\rbrace 
\leq \mathbb{P}\lbrace A \not \subset \phi \mid A_{1} \not \subset \phi \rbrace
\end{equation}
holds for all disjoint subsets $A$, $A_{1}$, $A_{2}$,   of $\mathcal{N}$
 such that  $\mathbb{P}\lbrace  A_{1} \not \subset \phi, A_{2} \not \subset \phi\rbrace  > 0 $. 
\end{theorem}
To get a similar result for two disjoint sets $A_{1}$, $A_{2}$ of $\mathcal{N}$, that is to prove that
 For $i_{0}, j_{0} \notin A_{1}\cup A_{2}$ we have
\begin{equation}\label{I-19}
\begin{split}
& P\{ i_{0}\in \phi \vert A_{1} \not \subset \phi, A_{2} \not \subset \phi\}
 \times 
P\{ j_{0} \in \phi \vert A_{1} \not \subset \phi, A_{2} \not \subset \phi \}\\
& - P\{ i_{0}\in \phi, j_{0} \in \phi \vert A_{1} \not \subset \phi, A_{2} \not \subset \phi \} \geq 0
\end{split}
\end{equation}
 is a more difficult task. Our proof rely upon a geometric/algebraic result  of an independent interest. To explain this, consider the subsets 
$\{a_{i}, i\in A_{1}\}$ and $\{a_{i}, i\in A_{2}\}$ of $\mathbb{R}^{p}$ where 
$ \mathfrak{A}=\{a^{1},...,a^{p}\} $, $ 1<p<N $, $1< p < N$, is a set of orthonormal vectors associated to the determinantal proces $\phi$. Without loss of generality we may suppose that:  $\vert A_{l}\vert <p$, 
$0< P\{A_{l}  \subset \phi\}=
\parallel \bigwedge_{i\in A_{l}}a_{i}\parallel ^{2}<1$, $l=1,2$ and obviously 
$$0<P\{ A_{1} \not \subset \phi, A_{2} \not \subset \phi\} 
= 1 - \parallel \bigwedge_{i\in A_{1}}a_{i}\parallel ^{2}
- \parallel \bigwedge_{i\in A_{2}}a_{i}\parallel ^{2}
+\parallel \bigwedge_{i\in A_{1}\cup A_{2}}a_{i}\parallel ^{2}<1.$$
With these properties in mind, we shall describe the algebraic result of interest which requires a thorough explanations.\\ 
Let $E_{l}$, $\vert E_{l}\vert=k_{l}<p$, $l=1,2$  be two subsets of $\mathbb{R}^{p}$  
 who fullfil the conditions above, that is:
 \begin{equation}\label{I-20}
 0<\kappa_{l}=\parallel \bigwedge_{x\in E_{l}}x\parallel ^{2}<1, l=1,2 
 \end{equation}
and
\begin{equation}\label{I-21}
0<1 - \parallel \bigwedge_{x\in E_{1}}x\parallel ^{2}
- \parallel \bigwedge_{i\in E_{2}}x\parallel ^{2}
+\parallel \bigwedge_{i\in E_{1}\cup E_{2}}x\parallel ^{2}<1.
\end{equation} 
Conditions (\ref{I-20}) implie that the subspaces  $\mathcal{E}_{l} $
 spanned by the subsets $E_{l}$, $l=1,2$ are respectively of dimension $k_{l}$.
Suppose $k_{1}\leq k_{2}$. There are two possible alternatives:\\
Case I:
 \begin{equation}\label{I-22}
\parallel \bigwedge_{i\in E_{1}\cup E_{2}}x\parallel ^{2}=0. 
\end{equation} 
In this case we have $\mathcal{E}_{1} \cap \mathcal{E}_{2}\neq \{0\}$ and it is
well-known \cite{ref1},  \cite{ref13}, \cite{ref22} that
there exist a set of orthonormal vectors 
$$u_{11},\dots,u_{1k}, u_{21},\dots, u_{2k},u_{31},\dots,u_{3k_{3}},u_{41},\dots,u_{4k_{4}},u_{51},\dots,u_{5k_{5}}$$  
(which is a basis of $\mathbb{R}^{p}$) with 
$$k+k_{3}=k_{1},\quad k+k_{3} + k_{4}=k_{2},\quad 
2k+k_{3}++k_{4}+k_{5}=p,$$  
as well as Jordan
principal angles $0<\alpha_{i}< \pi/2$, $i=1,\dots,k$\\
so that:\\
\begin{itemize}
\item[(i)]
the vectors $u_{11},\dots,u_{1k},u_{31},\dots, u_{3k_{3}}$ 
are a basis of $\mathcal{E}_{1}$
\item[(ii)]
Puting $v_{i}=u_{1i}\cos\alpha_{i} + u_{2i}\sin\alpha_{i}$
, $i=1,\dots,k$, the vectors 
$$v_{1},\dots,v_{k},u_{31},\dots,u_{3k_{3}},u_{41},\dots,u_{4k_{4}}
$$
are a basis of $\mathcal{E}_{2} $. 
\end{itemize}
This is the most general situation subject to the condition (\ref{I-22}). If $k_{1}=k_{2}$ 
then (obviously) the sequence $u_{4i}$, $i=1,\dots,k_{4}$ above must be taken out and if $2k+k_{3}+k_{4}=p$ then the same applies to the sequences
$u_{5i}$, $i=1,\dots,k_{5}.$\\
Case II: 
\begin{equation}\label{I-23}
\parallel \bigwedge_{i\in E_{1}\cup E_{2}}x\parallel ^{2}\neq 0. 
\end{equation} 
In this case we have $\mathcal{E}_{1} \cap \mathcal{E}_{2}= \{0\}$ and consequently the description above  works provided the sequence
$u_{3i}$, $i=1,\dots,k_{3}$ is taken out.
\\
Consider now two vectors $y,y'\in \mathbb{R}^{p}$, $y\neq y'$.
There exists 
$a=(a_{1},\dots,a_{k})$,   $b=(b_{1},\dots,b_{k})$, $c=(c_{1},\dots,c_{k_{3}})$, $d=(d_{1},\dots,d_{k_{4}})$
and $e=(e_{1},\dots,e_{k_{5})}$ such that writing
\begin{equation}\label{I-24}
\begin{split} 
 & y_{0}=\sum_{i=1}^{k}(a_{i}+b_{i}\cos\alpha_{i})u_{1i}
+ \sum_{i=1}^{k}b_{i}(\sin\alpha_{i})u_{2i}\\
 & y_{0}'=\sum_{i=1}^{k
 }(a_{i}'+b_{i}'\cos\alpha_{i})u_{1i}
+ \sum_{i=1}^{k}b_{i}'(\sin\alpha_{i})u_{2i}
\end{split}
\end{equation} 
we have
\begin{itemize}
\item[Case I.]
\begin{equation}\label{I-25}
\begin{split}
& y=y_{0}+ \sum_{i=1}^{k_{3}}c_{i}u_{3i}   
 + \sum_{i=1}^{k_{4}}d_{i}u_{4i}
 + \sum_{i=1}^{k_{5}}e_{i}u_{5i},\\
&  y'=y_{0}'+ \sum_{i=1}^{k_{3}}c_{i}'u_{3i}   
 + \sum_{i=1}^{k_{4}}d_{i}'u_{4i}
 + \sum_{i=1}^{k_{5}}e_{i}'u_{5i}
\end{split} 
\end{equation}
or 
\begin{equation}\label{I-26}
y=(y_{0},c,d,e),\quad y'=(y_{0}',c',d',e')
\end{equation}
in ccordinate representation for the ordered basis \\
 $u_{11},\dots,u_{1k}, u_{21},\dots, u_{2k},u_{31},\dots,u_{3k_{3}},u_{41},\dots,u_{4k_{4}},u_{51},\dots,u_{5k_{5}.}$ 
\item[Case II.]
\begin{equation}\label{I-27}
\begin{split}
& y=y_{0}  
 + \sum_{i=1}^{k_{4}}d_{i}u_{4i}
 + \sum_{i=1}^{k_{5}}e_{i}u_{5i},\\
&  y'=y_{0}' 
 + \sum_{i=1}^{k_{4}}d_{i}'u_{4i}
 + \sum_{i=1}^{k_{5}}e_{i}'u_{5i}
\end{split} 
\end{equation}
or 
\begin{equation}\label{I-28}
y=(y_{0},d,e),\quad y'=(y_{0}',d',e')
\end{equation}
in ccordinate representation for the ordered basis \\
 $u_{11},\dots,u_{1k}, u_{21},\dots, u_{2k},u_{41},\dots,u_{4k_{4}},u_{51},\dots,u_{5k_{5}}.$ 
\end{itemize}
Consider now
\begin{enumerate}
\item[Case I.]
 \begin{equation}\label{I-29}
\begin{split}
 \langle y,y'\rangle & - 
\kappa_{1}\langle y\wedge\bigwedge_{i=1}^{k}u_{1i}\bigwedge_{i=1}^{k_{3}}u_{3i},
y'\wedge\bigwedge_{i=1}^{k}u_{1i}\bigwedge_{i=1}^{k_{3}}u_{3i}\rangle\\
& - \kappa_{2}\langle y\wedge\bigwedge_{i=1}^{k}v_{i}\bigwedge_{i=1}^{k_{3}}u_{3i}\bigwedge_{i=1}^{k_{4}}u_{4i},  y'\wedge\bigwedge_{i=1}^{k}v_{i}\bigwedge_{i=1}^{k_{3}}u_{3i}\bigwedge_{i=1}^{k_{4}}u_{4i}\rangle.
\end{split} 
\end{equation}
\item[Case II.]
 \begin{equation}\label{I-30}
\begin{split}
 \langle y,y'\rangle & - 
\kappa_{1}\langle y\wedge\bigwedge_{i=1}^{k}u_{1i} ,
y'\wedge\bigwedge_{i=1}^{k}u_{1i}\rangle\\
& - \kappa_{2}\langle y\wedge\bigwedge_{i=1}^{k}v_{i} \bigwedge_{i=1}^{k_{4}}u_{4i},  y'\wedge\bigwedge_{i=1}^{k}v_{i}
\bigwedge_{i=1}^{k_{4}}u_{4i}\rangle\\
& + \kappa_{1}\kappa_{2}\prod_{i=1}^{k}(\sin\alpha_{i})^{2} \langle y\wedge\bigwedge_{i=1}^{k}u_{1i}\bigwedge_{i=1}^{k}v_{i}\bigwedge_{i=1}^{k_{4}}u_{4i},
y'\wedge\bigwedge_{i=1}^{k}u_{1i}\bigwedge_{i=1}^{k}v_{i}\bigwedge_{i=1}^{k_{4}}u_{4i}\rangle.
\end{split} 
\end{equation}
\end{enumerate}
The key point is the fact  - as previously seen in the special cases (\ref{I-12}) and (\ref{I-18}) -   that the expressions (\ref{I-29}) and (\ref{I-30} )  are  scalar products of two vectors $\tilde{y},\tilde{y}' \in \mathbb{R}^{p}$ : lemmas \ref{I-3}, \ref{I-8}, \ref{I-10} in section \ref{S4}. Furthermore  we will prove in Section \ref{S4} that:
\begin{theorem}\label{th2}
Define
\begin{enumerate}
\item[Case I.]
\begin{equation}\label{I-31}
\begin{split}
& \Lambda_{g1}(y,y')= \parallel\tilde{y}\wedge\tilde{y}' \parallel^{2}\\
 & - (1  - \kappa_{1} -\kappa_{2})
 \times
 \bigg[  \parallel y\wedge y'  \parallel^{2} -
 \kappa_{1}\parallel y\wedge y' \wedge \bigwedge_{i=1}^{k}u_{1i}\bigwedge_{i=1}^{k_{3}}u_{3i}\parallel^{2} \\ 
& \qquad -  \kappa_{2}\parallel y\wedge y' \wedge \ \bigwedge_{i=1}^{k}v_{i}\bigwedge_{i=1}^{k_{3}}u_{3i}\bigwedge_{i=1}^{k_{4}}u_{4i}
\parallel^{2}
 \bigg]. 
\end{split}
\end{equation}
\item[Case II.]
\begin{equation}\label{I-32}
\begin{split}
& \Lambda_{g2}(y,y')= \parallel\tilde{y}\wedge\tilde{y}' \parallel^{2}\\
 & - (1  - \kappa_{1} -\kappa_{2} + \kappa_{1}\kappa_{2}\prod_{i=1}^{k}(\sin\alpha_{i})^{2})
 \times
 \bigg[  \parallel y\wedge y'  \parallel^{2} -
 \kappa_{1}\parallel y\wedge y' \wedge \bigwedge_{i=1}^{k}u_{1i} \parallel^{2} 
 \\ 
& \qquad -  \kappa_{2}\parallel y\wedge y' \wedge \ \bigwedge_{i=1}^{k}v_{i} \bigwedge_{i=1}^{k_{4}}u_{4i}
\parallel^{2}\\
& \qquad + \kappa_{1}\kappa_{2}\prod_{i=1}^{k}(\sin\alpha_{i})^{2}
 \parallel y\wedge y' \wedge \ \bigwedge_{i=1}^{k}u_{1i}\bigwedge_{i=1}^{k}v_{i} \bigwedge_{i=1}^{k_{4}}u_{4i}
\parallel^{2}
 \bigg]. 
\end{split}
\end{equation}
\end{enumerate}
then
\begin{equation}\label{I-33}
\Lambda_{g1}(y,y')\geq 0 \quad  and \quad \Lambda_{g2}(y,y')\geq 0.
\end{equation}
 \end{theorem} 
 
This result provides a proof of the inequality (\ref{I-19}) from which, by Proposition \ref{t2} the conditional inequality (\ref{I-8}) for $n=3$ follows. Indeed, with the choice
$y=a_{i_{0}}$, $y'=a_{j_{0}}$, $E_{1}=\{a_{i}, i\in A_{1}\}$,
$E_{2}=\{a_{i}, i\in A_{2}\}$ we obtain:
\begin{itemize}
\item[(i)]
\begin{equation}\label{I-34}
\parallel \tilde{y} \parallel^{2}=P\{ i_{0}\in \phi, A_{1} \not \subset \phi, A_{2} \not \subset \phi\}, \quad
\parallel \tilde{y}' \parallel^{2}=P\{ j_{0}\in \phi, A_{1} \not \subset \phi, A_{2} \not \subset \phi\}
\end{equation}
\item[(ii)]
\begin{equation}\label{I-35}
\parallel \tilde{y} \wedge \tilde{y}'  \parallel^{2}=
P\{ i_{0}\in \phi, j_{0} \in \phi, A_{1} \not \subset \phi, A_{2} \not \subset \phi \}
\end{equation}
\item[(iii)]
Case I:
\begin{equation}\label{I-36}
P\{  A_{1} \not \subset \phi, A_{2} \not \subset \phi \} =
1  - \kappa_{1} -\kappa_{2} 
\end{equation}
\item[(iii')]
Case II:
\begin{equation}\label{I-37}
P\{  A_{1} \not \subset \phi, A_{2} \not \subset \phi \} =
1  - \kappa_{1} -\kappa_{2} + \kappa_{1}\kappa_{2}\prod_{i=1}^{k} \sin ^{2}\alpha_{i} 
\end{equation}
\end{itemize}
Thus, in virtue of (\ref{I-33}), substituting these probabilistic expressions into (\ref{I-31}) and (\ref{I-32})
and applying
the well-known formula 
\begin{equation}\label{I-38}
\parallel \tilde{y} \wedge \tilde{y}'  \parallel^{2}
=\parallel \tilde{y} \parallel^{2}\times \parallel  \tilde{y}'  \parallel^{2}
- \langle \tilde{y},\tilde{y}'\rangle^{2}
\end{equation}
 we obtain:
\begin{theorem}\label{th3}
For disjoint  subsets $\{i_{0}\}$, $\{j_{0}\}$, $A_{1}$ and $A_{2}$ of $\mathcal{N}$
such that  $\mathbb{P}\lbrace  A_{i} \not \subset \phi, i=1,2 \rbrace  > 0 $
we have: 
\begin{equation}\label{I-39}
\begin{split}
& P\{ i_{0}\in \phi \vert A_{1} \not \subset \phi, A_{2} \not \subset \phi\}
 \times 
P\{ j_{0} \in \phi \vert A_{1} \not \subset \phi, A_{2} \not \subset \phi \}\\
& - P\{ i_{0}\in \phi, j_{0} \in \phi \vert A_{1} \not \subset \phi, A_{2} \not \subset \phi \} \\
& =  P\{  A_{1} \not \subset \phi, A_{2} \not \subset \phi \}^{-2}\times \big( \langle \tilde{y},\tilde{y}'\rangle^{2} + 
\Lambda_{ig}(y,y')\big)\geq 0,i=1,2
\end{split}
\end{equation}
where  $P\{  A_{1} \not \subset \phi, A_{2} \not \subset \phi \}$ is given by (\ref{I-36}) and (\ref{I-37}). As a consequence the inequality
\begin{equation}
 \mathbb{P}\lbrace  A \not \subset \phi \mid A_{i} \not \subset \phi, \forall i=1,2,3\rbrace 
 \leq \mathbb{P}\lbrace A \not \subset \phi \mid A_{i} \not \subset \phi, \forall i=1,2\rbrace
\end{equation}
holds for all disjoint subsets $A$, $A_{i}$, $i=1,2,3$,   of $\mathcal{N}$
 such that  $\mathbb{P}\lbrace  A_{i} \not \subset \phi, \forall i=1,2, 3\rbrace  > 0 $. 
\end{theorem}
The rest of the paper is organized as follows. In section \ref{s2} we show
that in order to obtain (\ref{I-10}) it suffices to consider the case when the sets $A$ and $B$ are simple points. The idea of proof is pretty standart (see for example \cite{ref26}). It proceeds by induction and uses the fact that for a basic determinantal point process $\psi$ the conditionned process
$(\psi\vert y\in \phi)\backslash \lbrace y\rbrace$ is still a basic determinantal process. In section \ref{s3} we give a proof of the Theorem \ref{th1}. 
The section \ref{S4} is devoted to the proof of the main  Theorem \ref{th2} and consits of two stages. First, we deal in the subsection \ref{SS1}, with the case when $k=\vert E_{1}\vert =\vert E_{2}\vert$ and $p=2k$. This is the crucial
non-trivial step. In this particulat context we explicit the vectors  $\tilde{y}$, and   $\tilde{y}'$ whose scalar product coincides with the expression (\ref{I-30}) and obtain exact formulas for $\Lambda_{g2}$ which allows us to show that $\Lambda_{g2}\geq 0$. These results are extended  to the general cases (I and II) in the next subsection where exact formulas for $\tilde{y}$, and   $\tilde{y}'$, $\Lambda_{g1}$ and $\Lambda_{g2}$ are obtained as well: 
lemmas \ref{L8} - \ref{L11}. This  requires a "little  arithmetic" but does not offer particular difficulty. 
\section{ The reduction principle}\label{s2}
\begin{notation}
For disjoints sets $A_{i},i=1\dots n, $ such that\\ 
$P\{ A_{i} \not \subset \phi,i=1\dots n \}>0 $
we will note
\begin{equation}\label{40}
(\phi\vert A_{i} \not \subset \phi,i=1\dots n)= \psi.
\end{equation}
\end{notation}
\begin{prop}[reduction principle]\label{t2} 
Suppose that the inequality
\begin{equation}\label{41}
 P\{ x\in \psi, x' \in \psi   \}
 \leq P\{ x\in \psi  \}
 \times 
P\{ x' \in \psi \}
 \end{equation}
is satisfied for every choice of disjoint sets  $\lbrace x \rbrace$, $\lbrace x'\rbrace $, $A_{i},i=1\dots n $
and every basic discrete  determinantal process $\phi$, associated to a subspace of $ \mathbb{R}^{N} $, such that
$P\{ A_{i} \not \subset \phi,i=1\dots n \}>0 $.
Then we have:
\begin{equation}\label{42}
 P\{ B_{1}\subset \psi, B_{2} \subset \psi  \}
\leq P\{ B_{1}\subset \psi   \}
\times P\{ B_{2}\subset \psi  \} 
\end{equation}
for every choice of disjoint sets  $B_{1}$, $B_{2}$, $A_{i},i=1\dots n $
and every discrete  determinantal process $\phi$ such that
$P\{ A_{i} \not \subset \phi,i=1\dots n \}>0 $.
\end{prop}
Proof.-- By induction and the following lemma.
\begin{lemma}
Fix $1\leq k$ and disjoint sets  $B_{2} $, $A_{i} $, $i=1,\dots ,n$. Suppose (the induction hypothesis) that the inequality   (\ref{42}) holds for every choice of the set  $ B_{1}$ disjoint from  $B_{2}\cup(\bigcup_{i=1}^{n}A_{i}) $,   
and satisfying $\vert B_{1}\vert \leq k$. Then,  (\ref{42}) holds equally   under the condition $\vert B_{1}\vert \leq k+1 $.
\end{lemma}
Proof.-- Suppose that $\vert B_{1}\vert = k$ and let 
 $y\notin B_{1}\bigcup B_{2}\bigcup_{i=1}^{n}A_{i}$. We want to show that
\begin{equation}\label{43}
P\{ B_{1}\cup \lbrace y\rbrace\subset \psi, B_{2} \subset \psi   \}
\leq P\{ B_{1}\cup \lbrace y\rbrace \subset \psi   \}
\times P\{ B_{2}\subset \psi  \} .
\end{equation} 
We may suppose that $P\{ B_{1}\cup \lbrace y\rbrace\subset \phi, B_{2} \subset \phi , A_{i} \not \subset \phi,i=1\dots n \}>0$ (otherwise there in nothing to prove) which implies that $P\{y\in \phi\}>0$ and $P\{y\in \phi, A_{i} \not \subset \phi,i=1\dots n \}>0$. Furthermore, it is well known that 
\begin{equation}\label{44} 
\phi_{y} =(\phi\vert y\in \phi)\backslash \lbrace y\rbrace
\end{equation}
is a basic determinantal process (see, for example, property 3 of the proposition 1 in   \cite{ref12})

Consequently we have 
$$P\{ A_{i} \not \subset \phi_{y},i=1\dots n \}>0.$$
Induction hypothesis implies that
\begin{equation}\label{45}
\begin{split}
&P\{ B_{1}\subset \phi_{y}, B_{2} \subset \phi_{y} \vert A_{i} \not \subset \phi_{y},i=1\dots n \}
\\
& \leq P\{ B_{1}\subset \phi_{y} \vert A_{i} \not \subset \phi_{y},i=1\dots n \}
\times P\{ B_{2}\subset \phi_{y} \vert A_{i} \not \subset \phi_{y},i=1\dots n \} 
\end{split}
\end{equation}
which  by (\ref{44}) can be reformulated as follows:
\begin{equation}\label{46}
\begin{split}
& P\{ y \in \phi, B_{1}\subset \phi, B_{2} \subset \phi, A_{i} \not \subset \phi,i=1\dots n \}
\\
& \leq P\{y \in \phi, B_{1}\subset \phi,  A_{i} \not \subset \phi,i=1\dots n \}
\times
\dfrac{P\{y \in \phi, B_{2}\subset \phi,  A_{i} \not \subset \phi,i=1\dots n \}}{P\{y\in \phi, A_{i} \not \subset \phi,i=1\dots n \}} 
\end{split}
\end{equation}
$\Longleftrightarrow$
\begin{equation}\label{47}
P\{ y \in \psi, B_{1}\subset \psi, B_{2} \subset \psi  \}
\leq P\{y \in \psi, B_{1}\subset \psi  \}
\times
\dfrac{P\{y \in \psi, B_{2}\subset \psi  \}}{P\{y\in \psi \}}.
\end{equation}
From the induction hypothesis we get also
\begin{equation}\label{48}
P\{y \in \psi, B_{2}\subset \psi  \}\leq P\{y\in \psi \}\times
P\{ B_{2}\subset \psi  \}.
\end{equation}
This and   formulas (\ref{47}) provides (\ref{43}) as desired.
\section{Proof of theorem \ref{th1}}\label{s3}
\begin{itemize}
\item[1]
For $\vert A\vert = p$ such that  $0<\kappa =P\{ A\subset \phi\}<1$
we have trivially $$P\{ (\mathcal{N}\setminus A)\cap \phi\ \neq \emptyset \quad and \quad A \subset \phi\}=0.$$ 
Consequently,
\begin{equation}\label{49} 
\begin{split}
& P\{ i\in \phi \vert A \not \subset \phi\}
\times 
P\{ j \in \phi \vert A \not \subset \phi \} - P\{ i\in \phi, j \in \phi \vert A \not \subset \phi \} 
\\
& =
\dfrac{P\{ i\in \phi \}\times P\{ j\in \phi \}  
- P\{ A \not \subset \phi\}P\{ i\in \phi, j\in \phi\}}{P\{ A \not \subset \phi\}^{2
}}
\\
& 
= \dfrac{\parallel a_{i}\parallel^{2}\times \parallel a_{j}\parallel^{2}
- P\{ A \not \subset \phi\}\times \parallel a_{i}\wedge a_{j}\parallel^{2}}{(1-\kappa)^{2}}
\\
&
= \dfrac{ \langle a_{i},a_{j}\rangle^{2} + \kappa  \parallel a_{i}\wedge a_{j}\parallel^{2} }{ (1-\kappa)^{2}}
\end{split}
\end{equation}
where the last equality follows from the  formula 
$$\parallel a_{i}\wedge a_{j}\parallel^{2}=
\parallel a_{i}\parallel^{2}\times \parallel a_{j}\parallel^{2}
- \langle a_{i},a_{j}\rangle^{2}.$$
\item[2]
For $A=\{i_{1},\dots,i_{n-2}\} \subset \mathcal{N}$, $\vert A\vert=n-2 < p$,  $0<\kappa=P\{ A\subset \phi\}<1$,
choose a basis as  described in section \ref{S1} after the remark \ref{R2}. Then, with corresponding notations, namely  (\ref{I-13}), (\ref{I-14}) and (\ref{I-15})
we obtain
\begin{equation}\label{50}
\begin{split}
P\{ i\in \phi, A \not \subset \phi \} & =  P\{ i\in \phi \}
- P\{ i\in \phi, A  \subset \phi \}\\
& =
\parallel a_{i}\parallel^{2} - 
\parallel  a_{i}\wedge \bigwedge_{k=1}^{n-2}a_{i_{k}}\parallel^{2}
\\
& = \parallel\tilde{v}_{1}\parallel^{2} + \parallel W_{1}\parallel^{2}-
\parallel \bigwedge_{k=1}^{n-2}\cos\theta_{k}u^{k}\parallel^{2}\times \parallel W_{1}\parallel^{2}
\\
& = \parallel\tilde{v}_{1}\parallel^{2} + 
\parallel W_{1}\parallel^{2}\times (1- \kappa )
\end{split}
\end{equation}
 and similarly
\begin{equation}\label{51}
P\{ j\in \phi, A \not \subset \phi \}=\parallel\tilde{v}_{2}\parallel^{2} + 
\parallel W_{2}\parallel^{2}\times (1- \kappa ).
\end{equation}
Moreover, we have
\begin{equation}\label{52}
\begin{split}
P&\{ i\in \phi, j\in \phi, A \not \subset \phi \}  =
P\{ i\in \phi, j\in \phi \} -
P\{ i\in \phi, j\in \phi, A \subset \phi \} \\
& =
\parallel a_{i}\wedge a_{j}\parallel^{2} - 
\parallel a_{i}\wedge a_{j}\wedge \bigwedge_{k=1}^{n_2}a_{i_{k}}\parallel^{2}\\
& =\parallel a_{i}\wedge a_{j}\parallel^{2} - 
\kappa\parallel W_{1}\wedge W_{2} \parallel^{2}
 = \parallel a_{i}\parallel^{2}\times \parallel a_{j}\parallel^{2} 
- \langle a_{i},a_{j}\rangle^{2}\\
& \qquad -\kappa\big( \parallel W_{1}\parallel^{2}\times \parallel W_{2}\parallel^{2} - \langle W_{1},W_{2}\rangle^{2}\big) .
\end{split}
\end{equation}
Now, recall that $\parallel a_{i}\parallel^{2}= \parallel \tilde{v}_{1}\parallel^{2} 
+ \parallel W_{1}\parallel^{2}$, \\
$\parallel a_{j}\parallel^{2}= \parallel \tilde{v}_{2}\parallel^{2} 
+ \parallel W_{2}\parallel^{2}$ and $\langle a_{i},a_{j}\rangle =
\langle \tilde{v}_{1},\tilde{v}_{2}\rangle + \langle W_{1},W_{2}\rangle$.
From this and (\ref{50})-(\ref{52}) a simple calculation gives formula (\ref{I-17}). Besides, the obvious equality 
$\langle a_{i}\wedge(\bigwedge_{k=1}^{n-2} a_{i_{k}}), a_{j}\wedge(\bigwedge_{k=1}^{n-2} a_{i_{k})}\rangle
= \kappa \langle W_{1},W_{2} \rangle $ implies (\ref{I-18}).
\end{itemize} 
\section{Proof of Theorem \ref{th2}.}\label{S4}
\subsection{The case $k=\vert E_{1}\vert =\vert E_{2}\vert$
and $p=2k$.}\label{SS1}
Recal that  $E_{l}$, $\vert E_{l}\vert=k_{l}<p$, $l=1,2$  are two subsets of $\mathbb{R}^{p}$  who fullfil the conditions (\ref{I-20}) and (\ref{I-21}).
In addition we suppose that $k=\vert E_{1}\vert =\vert E_{2}\vert$
and $p=2k$. This and (\ref{I-20}) implie that the subspaces  $\mathcal{E}_{i} $
respectively spanned by the subsets $E_{i}$, $i=1,2$ are of dimension $k$
and that there exist a set of orthonormal vectors 
$u_{1},\dots,u_{k}, w_{1},\dots, w_{k}$ in $\mathbb{R}^{p}$ as well as Jordan
principal angles $0<\alpha_{i}< \pi/2$, $i=1,\dots,k$ so that the vectors
$u_{1},\dots,u_{k}$ 
are a basis of $\mathcal{E}_{1}$; the vectors $v_{i}=u_{i}\cos\alpha_{i} + w_{i}\sin\alpha_{i}$
, $i=1,\dots,k$, are a basis of 
$\mathcal{E}_{2}$: here, to simplify the notations we write $u_{i}$, $i=1,\dots,k$ (resp. $w_{i}$,$i=1,\dots,k$) instead $u_{1i}$, $i=1,\dots,k$ (resp. $u_{2i}$,$i=1,\dots,k$)
Hence,  we have:

\begin{equation}\label{53}
\begin{split}
& \bigwedge_{x\in E_{1}}x = \pm \sqrt{\kappa_{1}}\bigwedge_{i=1}^{k}u_{i},
\quad
\bigwedge_{x\in E_{2}}x = \pm\sqrt{ \kappa_{2}}\bigwedge_{i=1}^{k}v_{i},\\
& \bigwedge_{x\in E_{1}\cup E_{2}}x=\pm\sqrt{\kappa_{1}\kappa_{2}}
\prod_{i=1}^{k}\sin\alpha_{i}(\bigwedge_{i=1}^{k}u_{i})\wedge(\bigwedge_{i=1}^{k}w_{i})
\end{split}
\end{equation}
and
\begin{equation}\label{54}
\begin{split}
& 1 - \parallel \bigwedge_{x\in E_{1}}x\parallel ^{2}
- \parallel \bigwedge_{i\in E_{2}}x\parallel ^{2}
+\parallel \bigwedge_{i\in E_{1}\cup E_{2}}x\parallel ^{2}\\
& =
1 - \kappa_{1}- \kappa_{2} + 
\kappa_{1}\kappa_{2}\prod_{i=1}^{k}\sin^{2}\alpha_{i}.
\end{split}
\end{equation}
Introduce the notations
\begin{equation}\label{55}
\delta_{1,i}=1-\kappa_{1}\sin^{2}\alpha_{i}
\quad and \quad
\delta_{2,i}=1-\kappa_{2}\sin^{2}\alpha_{i}, i=1,\dots,k,
\end{equation}
and observe that by (\ref{I-20})  and (\ref{I-21}) we have, for all $i=1,\dots,k$, $0<\delta_{1,i}<1$, $0<\delta_{2,i}<1$
and $0<1 - \kappa_{1}- \kappa_{2} + 
\kappa_{1}\kappa_{2}\sin^{2}\alpha_{i}$. Consequently, we can define the angles $0<\alpha_{i}'<\pi/2$, $i=1,\dots,k$, via the formulas
 \begin{equation}\label{56}
\cos\alpha_{i}'= \frac{\cos\alpha_{i}}{\sqrt{\delta_{1,i}\delta_{2,i}}}
\end{equation}
and
\begin{equation}\label{57}
 \sin\alpha_{i}'
= \frac{\sin\alpha_{i}\sqrt{1-\kappa_{1}-\kappa_{2}
+\kappa_{1}\kappa_{2}\sin^{2}\alpha_{i}}}{\sqrt{\delta_{1,i}\delta_{2,i}}}.
\end{equation}
 \\
Now,  fix two vectors $y$, $y'$ in $\mathbb{R}^{p}$. They are of
the form
\begin{equation}\label{58}
y=\sum_{i=1}^{k}a_{i}u_{i} + \sum_{i=1}^{k}b_{i}v_{i}
=\sum_{i=1}^{k}(a_{i}+b_{i}\cos\alpha_{i})u_{i} + 
\sum_{i=1}^{k}b_{i}\sin\alpha_{i}w_{i}
\end{equation}
\begin{equation}\label{59}
y'=\sum_{i=1}^{k}a_{i}'u_{i} + \sum_{i=1}^{k}b_{i}'v_{i}
=\sum_{i=1}^{k}(a_{i}'+b_{i}'\cos\alpha_{i})u_{i} + 
\sum_{i=1}^{k}b_{i}'\sin\alpha_{i}w_{i}
\end{equation}
or, in ccordinate representation (for the ordered basis 
$u_{1},\dots,u_{k},w_{1}\dots,w_{k}$)

\begin{equation}\label{60}
y=(a_{1},\dots,a_{k},b_{1},\dots,b_{k})H, \quad
y'=(a_{1}',\dots,a_{k}',b_{1}',\dots,b_{k}')H 
\end{equation}
with the matrix

\begin{equation}\label{61}
H =
\left (\begin{array}{c|c}
I_{k} & \mathbf{0}  \\
\hline
\Delta & \Delta_{22}  \\
\end{array}\right)
\end{equation}

where $I_{k}$ is the unit matrix of size $k$ and 

$
\Delta=
\begin{pmatrix}
    \cos\alpha_{1} &  \dots & 0 \\
     \vdots &  \ddots & \vdots \\
    0 &  \dots & \cos\alpha_{k} 
\end{pmatrix}
$,
$
\Delta_{22}=
\begin{pmatrix}
    \sin\alpha_{1} &  \dots & 0 \\
     \vdots &  \ddots & \vdots \\
    0 &  \dots & \sin\alpha_{k} 
\end{pmatrix}
$
are diagonal matrices.
 
\begin{lemma}\label{L3}
Introduce 
\begin{enumerate}
\item
The matrix
$$
H_{1} =
\left (\begin{array}{c|c}
I_{k} & \mathbf{0}  \\
\hline
\Delta' & \Delta_{22}'  \\
\end{array}\right)
$$

where $I_{k}$ is the unit matrix of size $k$ and 

$
\Delta'=
\begin{pmatrix}
    \cos\alpha_{1}' &  \dots & 0 \\
     \vdots &  \ddots & \vdots \\
    0 &  \dots & \cos\alpha_{k}' 
\end{pmatrix},
$
$
\Delta_{22}'=
\begin{pmatrix}
    \sin\alpha_{1}' &  \dots & 0 \\
     \vdots &  \ddots & \vdots \\
    0 &  \dots & \sin\alpha_{k}' 
\end{pmatrix}
$
are 
\\
diagonal matrices.  
\item

The diagonal matrix

$
D=
\begin{pmatrix}
 \sqrt{\delta_{2,1}} &  \dots & & & & 0\\
 \vdots & \ddots & & \\
    &  & \sqrt{\delta_{2,k}}  & \\
    & & &\sqrt{\delta_{1,1}} &\\
    & & & & \ddots & \vdots\\
    0 & & & &\dots & \sqrt{\delta_{1,k}}
\end{pmatrix}
$
\end{enumerate}
Then, the vectors
\begin{equation}\label{62}
\tilde{y}=(a_{1},\dots,a_{k},b_{1},\dots,b_{k})D H_{1}
\end{equation}

\begin{equation}\label{63}
\tilde{y}'=(a_{1}',\dots,a_{k}',b_{1}',\dots,b_{k}')D H_{1}
\end{equation} 
or, in terms of the basis  $u_{1},\dots,u_{k}, w_{1},\dots, w_{k}$  
\begin{equation}\label{64}
\begin{split}
& \tilde{y}= \sum_{i=1}^{k}(a_{i}\sqrt{\delta_{2i}}+ 
b_{i}\sqrt{\delta_{1i}}\cos\alpha_{i}')u_{i} + 
\sum_{i=1}^{k}b_{i}\sqrt{\delta_{1i}} \sin\alpha_{i}'w_{i}\\
&\tilde{y}'=  \sum_{i=1}^{k}(a_{i}'\sqrt{\delta_{2i}}+ 
b_{i}'\sqrt{\delta_{1i}}\cos\alpha_{i}')u_{i} + 
\sum_{i=1}^{k}b_{i}'
\sqrt{\delta_{1i}} \sin\alpha_{i}'w_{i}
\end{split}
\end{equation}
satisfy the equality
\begin{equation}\label{65}
\langle \tilde{y},\tilde{y}'\rangle
=\langle y,y'\rangle - 
\langle y\wedge \bigwedge_{x\in E_{1}}x,
y'\wedge \bigwedge_{x\in E_{1}}x\rangle
-
\langle y\wedge \bigwedge_{x\in E_{2}}x,
y'\wedge \bigwedge_{x\in E_{2}}x \rangle. 
\end{equation}
\end{lemma}
\begin{proof}
According (\ref{53}), (\ref{58}), (\ref{59}) and noting that the multivectors
$w_{i}\wedge \bigwedge_{i=1}^{k}u_{i}$, $i=1,\dots,k$, are orthonormal we get
\begin{equation}\label{66}
\begin{split}
\langle y\wedge \bigwedge_{x \in E_{1}}x,y'\wedge \bigwedge_{x \in E_{1}}x \rangle 
& =\kappa_{1}\langle y\wedge \bigwedge_{i=1}^{k}u_{i},
y'\wedge \bigwedge_{i=1}^{k}u_{i} \rangle \\
& = \kappa_{1} \langle (\sum_{i=1}^{k}b_{i}\sin\alpha_{i}w_{i})\wedge \bigwedge_{i=1}^{k}u_{i},(\sum_{i=1}^{k}b_{i}\sin\alpha_{i}w_{i})\wedge \bigwedge_{i=1}^{k}u_{i}\rangle \\
& = \kappa_{1}\sum_{i=1}^{k}
b_{i}b_{i}'\sin^{2} \alpha_{i}.
\end{split}
\end{equation}
By the same argument we obtain
\begin{equation}\label{67}
\langle y\wedge \bigwedge_{x \in E_{2}}x,y'\wedge \bigwedge_{x \in E_{2}}x \rangle 
 = \kappa_{2}\sum_{i=1}^{k}
a_{i}a_{i}'\sin^{2} \alpha_{i}.
\end{equation}
From   (\ref{58}), (\ref{59}), (\ref{66}) and (\ref{67})  we get  
\begin{equation}\label{68}
\begin{split}
& \langle y,y'\rangle - 
\langle y\wedge \bigwedge_{x\in E_{1}}x,
y'\wedge \bigwedge_{x\in E_{1}}x\rangle
-
\langle y\wedge \bigwedge_{x\in E_{2}}x,
y'\wedge \bigwedge_{x\in E_{2}}x \rangle \\
& = \sum_{i=1}^{k}\left[ a_{i}a_{i}' (1 - 
\kappa_{2}\sin^{2} \alpha_{i}) + 
b_{i}b_{i}'(1 - 
\kappa_{1}\sin^{2} \alpha_{i}) + (a_{i}b_{i}' +
b_{i}a_{i}' )\cos\alpha_{i}\right] 
\end{split}
\end{equation}
which is exactly what a direct computation derived from the formulas (\ref{64}) gives.
\end{proof} 
Our goal now is to compute 
\begin{equation}\label{69}
\begin{split}
& \Lambda(y,y') = \parallel \tilde{y} \wedge \tilde{y}'  \parallel^{2}
  -
(1  - \kappa_{1} -\kappa_{2} + \kappa_{1}\kappa_{2}\prod_{i=1}^{k} \sin ^{2}\alpha_{i} )
\\ 
 & \qquad \qquad \times
 \left[  \parallel y\wedge y'  \parallel^{2}
  -   
\parallel y\wedge y'\wedge \bigwedge_{x\in E_{1}} x\parallel^{2}
-   
\parallel y\wedge y'\wedge \bigwedge_{x\in E_{2}} x \parallel^{2}\right]
\\
\end{split}
\end{equation}
and prove that $\Lambda(y,y')\geq 0$. This requires several stages.\\
First,  consider  the exterior product  
\begin{equation}\label{70}
\begin{split}
 ( a_{1},\dots,a_{k}&,b_{1},\dots,b_{k})\wedge
 (a_{1}',\dots,a_{k}',b_{1}',\dots,b_{k}')
\\
 &=\tilde{t}=(t(1,2),\dots,t(p-1,p))
\end{split}
 \end{equation}
the components of $\tilde{t}$ being indexed in the lexicographic order.
It is (certainly) well-known that we have
\begin{equation}\label{71}
\begin{split}
 y\wedge y'
 & =((a_{1},\dots,a_{k},b_{1},\dots,b_{k}) H ) \wedge ((a_{1}',\dots,a_{k}',b_{1},\dots,b_{k}' H) 
\\
& =\tilde{t} \widetilde{ H } 
\end{split}
\end{equation}
where $\widetilde{ H }$ is the compound matrix of order 2.
Properties of compound matrices (see \cite{gan} p.19 ) imply that
\begin{equation}\label{72}
\parallel y \wedge y'\parallel^{2}
 =\parallel \tilde{t} \widetilde{ H } \parallel^{2}
= \langle \tilde{t},\tilde{t}  
 \widetilde{H H^{t} } \rangle.
\end{equation}
Straighforwardly
\begin{equation}\label{73}
H H^{t} =
\left (\begin{array}{c|c}
I_{k} & \Delta  \\
\hline
\Delta & I_{k} \\
\end{array}\right)
\end{equation}
and thus the compound matrix $\widetilde{H H^{t} } $ has a  simple structure.
Precisely; for $1\leq i <j\leq k$ denote
 \begin{equation}\label{74}
\overline{t}(i,j) = (t(i,j),t(i,k+j),t(j,k+i),t(k+i,k+j))
\end{equation}
and accordingly consider the following order 
\begin{equation}\label{75}
\begin{split}
& (1,2) (1,2+k) (2,1+k) (1+k,2+k)....(i,j)(i,j+k)(j,i+k)(i+k,j+k)
\\
& ...(k-1,k) (k-1,2k)(k,2k-1)(2k-1,2k)(1,1+k)(2,2+k).....(k,2k).
\end{split}
\end{equation}
obtained by concatenating the  $\overline{t}(i,j)$ $1\leq i<j\leq k$ in the
lexicographic order followed by the sequence  ($1,1+k)(2,2+k).....(k,2k)$.
Then we have:
\begin{lemma}\label{L4}. 
By reordering the row and columns of  $\widetilde{H H^{t} }$ according to the order given above  we obtain the following block-diagonal matrix
\begin{equation}\label{76}
\begin{pmatrix}
    A(1,2) &  \dots &   &  & &  \\
  \vdots    & \ddots  &  & & &\\
     &     &   A(k-1,k)  &  & &\\
   &      &    &  \sin^{2}\alpha_{1}& & \\
   &   &   &   & \ddots& \vdots  & \\
   &  &   & &     \dots   &  \sin^{2}\alpha_{k}\\
 \end{pmatrix}
\end{equation}
whith
\begin{equation}\label{77}
A(i,j)=
 \begin{pmatrix}
    1 &  \cos\alpha_{j} & -\cos\alpha_{i}  & \cos\alpha_{i}\cos\alpha_{j} \\
  \cos\alpha_{j}     &  1 & - \cos\alpha_{i}\cos\alpha_{j}&\cos\alpha_{i} \\
 -\cos\alpha_{i}    &  -\cos\alpha_{i}\cos\alpha_{j}  &  1 & -\cos\alpha_{j} \\
 \cos\alpha_{i}\cos\alpha_{j}   &   \cos\alpha_{i}      & -\cos\alpha_{j}   & 1\\
 \end{pmatrix}.
\end{equation}
\end{lemma}
\begin{proof}
With the notations of \cite{gan} denote by 
$(\begin{array}{cc}
         j_{1} & j_{2} \\
         i_{1} & i_{2} \\
  \end{array})
 $
the minors of order 2 of the matrix $H H^{t}$ and define for $1\leq i <j \leq k$
the matrices
$$
A(i,j)= \left[ 
\begin{array}{c@{}c@{}c@{}c}
 (\begin{array}{cc}
 i & j \\
         i & j \\
 \end{array})
  &  (\begin{array}{cc}
                        i & j \\ 
                       i & k+j \\
                      \end{array})
                      & (\begin{array}{cc}
         i & j \\
         j & k+i \\
  \end{array})
  &(\begin{array}{cc}
                        i & j \\ 
                       k+i & k+j \\
                      \end{array})                   
                      \\
  (\begin{array}{cc}
         i & k+j \\
         i & j \\
  \end{array}) & (\begin{array}{cc}
                        i & k+j \\ 
                       i & k+j \\
                      \end{array}) 
                      &
                      (\begin{array}{cc}
         i & k+j \\
         j & k+i \\
  \end{array})
  &(\begin{array}{cc}
           i & k+j \\
         k+i & k+j  \\
  \end{array})
   \\
   (\begin{array}{cc}
         j & k+i \\
         i & j \\ 
   \end{array})        
          & (\begin{array}{cc}
                       j & k+i \\ 
                       i & k+j \\
                  \end{array}) & (\begin{array}{cc}
         j & k+i \\
         j & k+i  \\
  \end{array}) & (\begin{array}{cc}
         j & k+i \\
         k+i & k+j  \\
  \end{array})
  \\
(\begin{array}{cc}
         k+i & k+j \\
         i & j \\
  \end{array}) & (\begin{array}{cc}
         k+i & k+j  \\
           i & k+j \\
  \end{array}) & \ (\begin{array}{cc}
                                   k+i & k+j \\
                                   j & k+i \\
                                  \end{array}) &(\begin{array}{cc}
         k+i & k+j \\
         k+i & k+j  \\
  \end{array})
                                  \\
\end{array} \right] 
$$
It is obvious that the  not null minors of the matrix $H H^{t}$ are given either by the elements of matrices $A(i,j)$, $1\leq i <j \leq k$
or are of the form 
$$
(\begin{array}{cc}
i&  k+i\\
i&  k+i\\
\end{array})
=
\det
\begin{pmatrix}
1& \cos\alpha_{i} \\
\cos\alpha_{i}& 1\\
\end{pmatrix}
=\sin^{2}\alpha_{i}.
$$
Therefore, we obtain (\ref{76}), and computing the elements of the matrices $A(i,j)$, $1\leq i <j \leq k$,  (\ref{77}) follows.
\end{proof}
By means of lemma \ref{L4} and (\ref{72}) we deduce:
\begin{lemma}\label{L5}
We have: 
\begin{equation}\label{78}
\begin{split}
& \parallel y\wedge y'\parallel^{2}
-\parallel y\wedge y'\wedge \bigwedge_{x\in E_{1}}x\parallel^{2}
-\parallel y\wedge y'\wedge \bigwedge_{x\in E_{2}}x\parallel^{2}
\\
& = \sum_{1\leq i<j\leq k}\langle \tilde{t}(i,j),\tilde{t}(i,j)B(i,j)\rangle
+ \sum_{1\leq i \leq k}t(i,k+i)^{2}\sin^{2}\alpha_{i}
\end{split}
\end{equation}
where the symmetric matrices
$
B(i,j)=
$
\begin{equation}\label{79}
 \begin{pmatrix}
    1 - \kappa_{2}(\sin\alpha_{i}\sin\alpha_{j})^{2} &  \cos\alpha_{j} & -\cos\alpha_{i}  & \cos\alpha_{i}\cos\alpha_{j} \\
  \cos\alpha_{j}     &  1 & -\cos\alpha_{i}\cos\alpha_{j}&\cos\alpha_{i} \\
 -\cos\alpha_{i}    &  -\cos\alpha_{i}\cos\alpha_{j}  &  1 & -\cos\alpha_{j} \\
 \cos\alpha_{i}\cos\alpha_{j}   &   \cos\alpha_{i}      & -\cos\alpha_{j}   & 1  - \kappa_{1}(\sin\alpha_{i}\sin\alpha_{j})^{2} \\
 \end{pmatrix}
\end{equation}
are  positive-definite.
\end{lemma}
\begin{proof}
From (\ref{72}) and (\ref{74}) - (\ref{76}) we get
\begin{equation}\label{80}
\begin{split}
\parallel y \wedge y'\parallel^{2}
 =\parallel \tilde{t} \widetilde{ H } \parallel^{2}
&  = \langle \tilde{t},\tilde{t}  
 \widetilde{H H^{t} } \rangle 
  \\
& = \sum_{1\leq i<j\leq k}\langle \tilde{t}(i,j),\tilde{t}(i,j)A(i,j)\rangle
+ \sum_{1\leq i \leq k}t(i,k+i)^{2}\sin^{2}\alpha_{i}.
\end{split}
\end{equation}
Moreover, from (\ref{53}) (\ref{58}) (\ref{59}) we obtain
\begin{equation}\label{81}
\begin{split}
 & y\wedge y'\wedge \bigwedge_{x\in E_{1}}x = 
 \pm \sqrt{\kappa_{1}}(\sum_{i=1}^{k}b_{i}
 v_{i})\wedge(\sum_{i=1}^{k}b_{i}'
 v_{i})\wedge \bigwedge_{i=1}^{k}u_{i}
 \\
 & = \pm \sqrt{\kappa_{1}}\left[ \sum_{i=1}^{k}b_{i}
 (u_{i}\cos \alpha_{i}+w_{i}\sin \alpha_{i})\right]
  \wedge \left[ \sum_{i=1}^{k}b_{i}'
 (u_{i}\cos \alpha_{i}+w_{i}\sin \alpha_{i})\right]
 \\
 & \qquad \wedge \bigwedge_{i=1}^{k}u_{i}
 \\
 & =
 \pm \sqrt{\kappa_{1}} \sum_{1\leq i < j \leq k}
 (b_{i}b_{j}' - b_{i}'b_{j})(\sin\alpha_{i}\sin\alpha_{j})^{2}
 w_{i}\wedge w_{j}
   \wedge \bigwedge_{i=1}^{k}u_{i}.
\end{split}
\end{equation}
the multivectors apearing above being orthonormal. Therefore
\begin{equation}\label{82}
 \parallel y\wedge y'\wedge \bigwedge_{x\in E_{1}}x \parallel^{2}
 = \kappa_{1}\sum_{1\leq i < j \leq k}t(k+i,k+j)^{2}(\sin\alpha_{i}\sin\alpha_{j})^{2}.
\end{equation}
By a similar computation we have also
\begin{equation}\label{83}
 \parallel y\wedge y'\wedge \bigwedge_{x\in E_{2}}x \parallel^{2}
 = \kappa_{2}\sum_{1\leq i < j \leq k}t(i,j)^{2}(\sin\alpha_{i}\sin\alpha_{j})^{2}
\end{equation}
Clearly, (\ref{80}), (\ref{82}) and (\ref{83}) implie (\ref{78}) and (\ref{79}).\\
In order to see that the matrices $B(i,j)$ are positive definite it suffices to show that their leading principal minors are positive. Recall that 
$0<\kappa_{1} <1$, $i=1,2$ and that $0<1-\kappa_{1}-\kappa_{2}
+\kappa_{1}\kappa_{2}\prod_{i=1}^{k}\sin^{2}\alpha_{i}<1$. Computing  we obtain
\begin{itemize}
\item[(i)]
\begin{equation}\label{84}
\begin{split}
 \det& B(i,j)
\\
& = 
(\sin\alpha_{i}\sin\alpha_{j})^{2}\left[
1- \kappa_{1} -\kappa_{2} +\kappa_{1}\kappa_{2}
(1 - (\cos\alpha_{i}\cos\alpha_{j})^{2})\right] >0
\end{split}
\end{equation}
\item[(ii)]
\begin{equation}\label{85}
\begin{split}
& \det
\begin{pmatrix}
    1 - \kappa_{2}(\sin\alpha_{i}\sin\alpha_{j})^{2} &  \cos\alpha_{j} & -\cos\alpha_{i}   \\
  \cos\alpha_{j}     &  1 & -\cos\alpha_{i}\cos\alpha_{j}  \\
 -\cos\alpha_{i}    &  -\cos\alpha_{i}\cos\alpha_{j}  &  1  \\
\end{pmatrix}
\\
& =
(\sin\alpha_{i}\sin\alpha_{j})^{2}
\left[
1- \kappa_{2}
(1 - (\cos\alpha_{i}\cos\alpha_{j})^{2} )\right]  >0
\end{split}
\end{equation}
\item[(iii)]
\begin{equation}\label{86}
\begin{split}
& \det
\begin{pmatrix}
1 - \kappa_{2}(\sin\alpha_{i}\sin\alpha_{j})^{2}     &   \cos\alpha_{j}   \\
  \cos\alpha_{j}    &  1  \\
\end{pmatrix}
\\
& =\sin^{2}\alpha_{j}( 1 - \kappa_{1}\sin^{2}\alpha_{i})  >0
\end{split}
\end{equation}
\item[(iv)]
\begin{equation}\label{87}
1 - \kappa_{2}(\sin\alpha_{i}\sin\alpha_{j})^{2} >0
\end{equation}
\end{itemize}
 as desired.
\end{proof}
Now we proceed to evaluate the norm of the exterior product $\tilde{y}\wedge \tilde{y}'$.
As above we have 
\begin{equation}\label{88}
\begin{split}
 \tilde{y}\wedge \tilde{y}'
 & =((a_{1},\dots,a_{k},b_{1},\dots,b_{k})D H_{1} ) \wedge ((a_{1}',\dots,a_{k_{1}}',b_{1}',\dots,b_{k}')D H_{1}) 
\\
& =\tilde{t} \widetilde{D H_{1} } 
\end{split}
\end{equation}
where $\widetilde{D H_{1} }=\widetilde{D}\widetilde{H_{1}
}$ is the compound matrix of order 2
and therefore
\begin{equation}\label{89}
\parallel \tilde{y} \wedge\tilde{y}'\parallel^{2}
 =\parallel \tilde{t} \widetilde{D H_{1} } \parallel^{2}
=\langle \tilde{t},\tilde{t} \tilde{D}
\widetilde{H_{1} H_{1}^{t} }\tilde{D}\rangle. 
\end{equation}
The descripton of the matrices $H_{1}$ and $D$ given in
Lemma \ref{L3} shows that the matrix $\tilde{D}
\widetilde{H_{1} H_{1}^{t}}\tilde{D}
=(DH_{1} H_{1}^{t}D\widetilde{)} $ 
has exactly the same structure as the matrix 
$\widetilde{HH^{t}}$. Hence we obtain trivially.
\begin{lemma}\label{L6}
Reordering the row and colums of the matrix $\tilde{D}
\widetilde{H_{1} H_{1}^{t}}\tilde{D}
$ as described in lemma \ref{L4} we obtain mutatis mutandis the 
following block-diagonal matrix
$$
\begin{pmatrix}
    A'(1,2) &  \dots &   &  & &  \\
  \vdots    & \ddots  &  & & &\\
     &     &   A'(k-1,k)  &  & &\\
   &      &    &  s(1)& & \\
   &   &   &   & \ddots& \vdots  & \\
   &  &   & &     \dots   &  s(k)\\
 \end{pmatrix}
 $$
with
\begin{equation}\label{90}
\begin{split}
&A'(i,j)=
\\
&
\begin{pmatrix}
    \delta_{2i}\delta_{2j} & \delta_{2i} \cos\alpha_{j} & -\delta_{2j}\cos\alpha_{i}  & \cos\alpha_{i}\cos\alpha_{j} \\
 \delta_{2i} \cos\alpha_{j}     &  \delta_{2i}\delta_{1i} & -\cos\alpha_{i}\cos\alpha_{j}&\delta_{1j}\cos\alpha_{i} \\
 -\delta_{2j}\cos\alpha_{i}    &  -\cos\alpha_{i}\cos\alpha_{j}  &  \delta_{2j}\delta_{1i} & -\delta_{1i}\cos\alpha_{j} \\
 \cos\alpha_{i}\cos\alpha_{j}   &   \delta_{1j}\cos\alpha_{i}      & -\delta_{1i}\cos\alpha_{j}   & \delta_{1i}\delta_{1j}\\
 \end{pmatrix}
 \end{split}
\end{equation}

 and
\begin{equation}\label{91}
 s(i)=(1- \kappa_{1}- \kappa_{2}+\kappa_{1}\kappa_{2}\sin^{2}\alpha_{i})\sin^{2}\alpha_{i}.
\end{equation}
Consequently
\begin{equation}\label{92}
\begin{split}
\parallel \tilde{y} \wedge\tilde{y}'\parallel^{2}
 & =\langle \tilde{t},\tilde{t} \tilde{D}
\widetilde{H_{1} H_{1}^{t} }\tilde{D}\rangle  
  \\
& = \sum_{1\leq i<j\leq k}\langle \tilde{t}(i,j),\tilde{t}(i,j)A'
(i,j))\rangle
\\
& + \sum_{1\leq i \leq k}t(i,k+i)^{2}(1- \kappa_{1}- \kappa_{2}+\kappa_{1}\kappa_{2}\sin^{2}\alpha_{i})\sin^{2}\alpha_{i}.
\end{split}
\end{equation}

\end{lemma}
\begin{remark}\label{R4}
 Notice that the elements of the matrix $A'(i,j)$, as well as the  $s(i)$,
 which are initially expressed in terms of angles $\alpha_{i}'$, have been
 recalculated by means of formulas (\ref{I-25})-(\ref{I-26}) in terms of  $\alpha_{i}$. 
\end{remark}

Now, we are in position  to compute the exact value of $\Lambda(y,y')$.  
\begin{theorem}\label{Th4}
Consider the matrices 
\begin{equation}\label{93}
M_{ij}=
 \begin{pmatrix}
  1  - \kappa_{2}(\sin\alpha_{i}\sin\alpha_{j})^{2}  & \cos\alpha_{i}\cos\alpha_{j} \\
 \cos\alpha_{j} & \cos\alpha_{i}& \\
 -\cos\alpha_{i} & -\cos\alpha_{j}& \\
 \cos\alpha_{i}\cos\alpha_{j}   &  1  - \kappa_{1}
 (\sin\alpha_{i}\sin \alpha_{j})^{2}\\
 \end{pmatrix}
\end{equation}
$1\leq i <j\leq k$ and the diagonal matrix
\begin{equation}\label{94}
\Delta_{\kappa}=
 \begin{pmatrix}
 \\ \sqrt{\kappa_{1}}    & 0 \\
 0     &   \sqrt{\kappa_{2}}   \\
\end{pmatrix}.
\end{equation}
Then 
\begin{enumerate}
\item[1]
For $ k= 2$ we have
\begin{equation}\label{95}
\begin{split}
\Lambda(y,y')= & \parallel \tilde{t}(1,2)M_{12}\Delta_{\kappa} \parallel ^{2}
+ \\
 & \kappa_{1}\kappa_{2}(\sin\alpha_{1}\sin\alpha_{2})^{2}\big[ t(1,3)
 \sin\alpha_{2} 
+ t(2,4)\sin\alpha_{1}\big]^{2}
 \end{split}
\end{equation}
\item[2.]
For $ k> 2$ we have
\begin{equation}\label{96}
\Lambda(y,y')=  \Lambda_{1}(y,y')+\Lambda_{2}(y,y')+\Lambda_{3}(y,y')
\end{equation}
where
\begin{enumerate}
\item[(i)]
\begin{equation}\label{97}
\Lambda_{1}(y,y')=\sum_{1\leq i<j\leq k}\parallel \tilde{t}(i,j)M_{ij}\Delta_{\kappa} \parallel ^{2}
\end{equation}
\item[(ii)]
\begin{equation}\label{98}
\begin{split}
 & \Lambda_{2}(y,y')= \\
& \kappa_{1}\kappa_{2}\bigg[ \sum_{1\leq i \leq k}t(i,k+i)^{2}
\bigg(1 - \prod_{1\leq j\leq k, j\neq i} \sin ^{2}\alpha_{j} \bigg)
\sin^{4}\alpha_{i}
\\
& +  2\sum_{1\leq i<j  \leq k}t(i,k+i)t(j,j+k)
(\sin\alpha_{i}\sin\alpha_{j})^{2}\cos\alpha_{i}\cos\alpha_{j}\bigg]
\end{split} 
\end{equation}
\item[(iii)]
\begin{equation}\label{99}
\begin{split}
&  \qquad \Lambda_{3}(y,y') = \\
& \kappa_{1}\kappa_{2}\sum_{1\leq i<j\leq k}
(\sin\alpha_{i}\sin\alpha_{j})^{2}
\bigg(1 - \prod_{1\leq l\leq k, l\neq i,j} \sin ^{2}\alpha_{l} \bigg) 
\langle \tilde{t}(i,j),\tilde{t}(i,j)B
(i,j)\rangle
\end{split} 
\end{equation}
\end{enumerate}
satisfy $\Lambda_{i}(y,y')\geq 0$, $i=1,2,3$.
 \end{enumerate}
\end{theorem}
\begin{proof}
The case $k=2$.\\
 From (\ref{78}) and (\ref{92})we obtain
\begin{equation}\label{100}
\begin{split}
& \Lambda(2)=
\langle \tilde{t}(1,2),\tilde{t}(1,2)A'
(1,2)\rangle
\\
& +
\sum_{i=1}^{2}(t(i,2+i))^{2}(1- \kappa_{1}- \kappa_{2}+\kappa_{1}
\kappa_{2}\sin^{2}\alpha_{i})\sin^{2}\alpha_{i}.
\\
& - 
(1  - \kappa_{1} -\kappa_{2} + \kappa_{1}\kappa_{2}(
\sin \alpha_{1} \sin \alpha_{2})^{2})\times \bigg[
\langle \tilde{t}(1,2),\tilde{t}(1,2)B(1,2)\rangle
\\
& 
 \quad \quad + \sum_{i=1}^{2}(t(i,2+i))^{2}\sin^{2}\alpha_{i}\bigg]
\\
& =\langle \tilde{t}(1,2),\tilde{t}(1,2)C(1,2)\rangle
 + \kappa_{1}\kappa_{2}\bigg[t(1,3)^{2}\sin^{4}\alpha_{1}
\cos^{2}\alpha_{2}
\\
&
 + t(2,4)^{2}\sin^{4}\alpha_{2}
\cos^{2}\alpha_{1}\bigg] 
\end{split}
\end{equation}
where
$$
C(1,2)=A'(1,2) -(1  - \kappa_{1} -\kappa_{2} + \kappa_{1}\kappa_{2} (\sin\alpha_{1}\sin\alpha_{2} )^{2})B(1,2).
$$
Now, computing from (\ref{79}) and  (\ref{90})   the elements of the  symmetric matrix 
$C(1,2)= (c^{i}_{j})_{i,j=1,\dots,4}
$ 
we obtain:
\begin{equation}\label{101}
\begin{split}
& c^{1}_{1}=\kappa_{2}(\cos\alpha_{1}\cos\alpha_{2})^{2}
+ \kappa_{1}\bigg(1 - \kappa_{2}(\sin\alpha_{1}\sin\alpha_{2})^{2}\bigg)^{2}\\
& c^{1}_{2}=\cos\alpha_{2}\left[ 
\kappa_{2}(\cos\alpha_{1})^{2} + \kappa_{1}(1 - \kappa_{2}(\sin\alpha_{1}\sin\alpha_{2})^{2})
\right]\\ 
& c^{1}_{3}=-\cos\alpha_{1}\left[ 
\kappa_{2}(\cos\alpha_{2})^{2} + \kappa_{1}(1 - \kappa_{2}(\sin\alpha_{1}\sin\alpha_{2})^{2})
\right]\\ 
& c^{1}_{4}=\cos\alpha_{1}\cos\alpha_{2}\left[
\kappa_{2}  + \kappa_{1}(1 - \kappa_{2}(\sin\alpha_{1}\sin\alpha_{2})^{2})
\right] \\
& c^{2}_{2}=\kappa_{1}(\cos\alpha_{2})^{2} + 
\kappa_{2}(\cos\alpha_{1})^{2}\\
& c^{2}_{3}=- \cos\alpha_{1}\cos\alpha_{2}\left[
\kappa_{2}  + \kappa_{1}(1 - \kappa_{2}(\sin\alpha_{1}\sin\alpha_{2})^{2})
\right]\\
& c^{2}_{4}=\cos\alpha_{1}
\left[ 
\kappa_{1}(\cos\alpha_{2})^{2} + \kappa_{2}(1 - \kappa_{1}(\sin\alpha_{1}\sin\alpha_{2})^{2})
\right]\\
& c^{3}_{3}=\kappa_{1}(\cos\alpha_{1})^{2} + 
\kappa_{2}(\cos\alpha_{2})^{2} \\
& c^{3}_{4}=-\cos\alpha_{2}
\left[ 
\kappa_{1}(\cos\alpha_{1})^{2} + \kappa_{2}(1 - \kappa_{1}(\sin\alpha_{1}\sin\alpha_{2})^{2})
\right]\\ 
& c^{4}_{4}=\kappa_{1}(\cos\alpha_{1}\cos\alpha_{2})^{2}
+ \kappa_{2}\bigg(1 - \kappa_{1}(\sin\alpha_{1}\sin\alpha_{2})^{2}\bigg)^{2} 
\end{split}
\end{equation} 

Recall (\ref{74}) that $\tilde{t}(1,2)=(t(1,2),t(1,4),t(2,3),t(3,4))$. A close look at  formulas above enables us to see that
\begin{equation}\label{102}
\begin{split}
&  \langle \tilde{t}(1,2),\tilde{t}(1,2)
C(1,2))\rangle 
  =
\kappa_{2}\Bigl[ t(1,2)\cos\alpha_{1}\cos\alpha_{2}
+ t(1,4)\cos\alpha_{1} - t(2,3)\cos\alpha_{2}
\\
& \qquad + t(3,4)(1 - \kappa_{1}(\sin\alpha_{1}\sin\alpha_{2})^{2})
\Bigr]^{2} + \kappa_{1}\Bigl[ t(1,2)(1 - \kappa_{2}(\sin\alpha_{1}\sin\alpha_{2})^{2})
\\
& \qquad + t(1,4)\cos\alpha_{2} - t(2,3)\cos\alpha_{1}
+ t(3,4)\cos\alpha_{1}\cos\alpha_{2}\Bigr]^{2}
\\
& \qquad + 2\kappa_{1}\kappa_{2}
(\sin\alpha_{1}\sin\alpha_{2})^{2}\cos\alpha_{1}\cos\alpha_{2}\Bigl( t(1,2)t(3,4)
+ t(1,4)t(2,3)
\Bigr).
\end{split} 
\end{equation}
Applying now the Pl\'ucker relation
\begin{equation}\label{103}
t(1,2)t(3,4)
+ t(1,4)t(2,3)=t(1,3)t(2,4).
\end{equation}
to the right side of (\ref{102}) we obtain 
from (\ref{100}), (\ref{102}) and (\ref{103}) that
\begin{equation}\label{104}
\begin{split}
&  \Lambda(2)=
 \kappa_{2}\Bigl[ t(1,2)\cos\alpha_{1}\cos\alpha_{2}
+ t(1,4)\cos\alpha_{1} - t(2,3)\cos\alpha_{2}
\\
& \qquad + t(3,4)(1 - \kappa_{1}(\sin\alpha_{1}\sin\alpha_{2})^{2})
\Bigr]^{2} + \kappa_{1}\Bigl[ t(1,2)(1 - \kappa_{2}(\sin\alpha_{1}\sin\alpha_{2})^{2})
\\
& \qquad + t(1,4)\cos\alpha_{2} - t(2,3)\cos\alpha_{1}
+ t(3,4)\cos\alpha_{1}\cos\alpha_{2}\Bigr]^{2}
\\
& \qquad + \kappa_{1}\kappa_{2}\Bigl[
t(1,3)\sin^{2}\alpha_{2}
\cos\alpha_{1} + t(2,4)\sin^{2}\alpha_{1}
\cos\alpha_{2}
\Bigr]^{2}\\
& =\parallel \tilde{t}(1,2)M_{12}\Delta_{k} \parallel ^{2}
+ \\
 & \qquad \kappa_{1}\kappa_{2}\big[ t(1,3)\sin^{2}\alpha_{2}\sin\alpha_{1}  
+ t(2,4)\sin^{2}\alpha_{1}\sin\alpha_{2}\big]^{2}
\end{split} 
\end{equation}
as desired.\\
Now we will proceed for k >2.\\
From (\ref{69}), (\ref{78}) and (\ref{92}) we obtain
\begin{equation}\label{105}
\begin{split}
 \Lambda(y,y')  = & \sum_{1\leq i<j \leq k}\langle \tilde{t}(i,j),\tilde{t}(i,j)A'
(i,j)\rangle
\\
 & + \sum_{1\leq i \leq k}t(i,k+i)^{2}(1- \kappa_{1}- \kappa_{2}+\kappa_{1}\kappa_{2}\sin^{2}\alpha_{i})\sin^{2}\alpha_{i} 
\\
&   -
(1  - \kappa_{1} -\kappa_{2} + \kappa_{1}\kappa_{2}\prod_{i=1}^{k} \sin ^{2}\alpha_{i} )
\times \bigg[ 
\sum_{1\leq i<j \leq k}\langle \tilde{t}(i,j),\tilde{t}(i,j)B(i,j)\rangle
\\
&   + \sum_{1\leq i \leq k}t(i,k+i)^{2}\sin^{2}\alpha_{i}
\bigg].
\end{split}
\end{equation}
Denote
$$
C(i,j)=A'(i,j)- (1- \kappa_{1} - \kappa_{2}
+  \kappa_{1}\kappa_{2}(\sin\alpha_{i}\sin\alpha_{j})^{2})B(i,j).
$$
Then, substituting in (\ref{101} ) the indices $i,j$ in place of the indices $1,2$ we obtain the elements of the matrix $C(i,j)$ and proceeding as for $C(1,2)$, 
using the Pl\'ucker relation $t(i,j)t(i+k,j+k) + t(i,j+k)t(j,i+k)=t(i,i+k)t(j,j+k) $,
this gives (with notation (\ref{97}))
\begin{equation}\label{106}
\begin{split}
& \sum_{1\leq i<j\leq k}\langle \tilde{t}(i,j),\tilde{t}(i,j)C
(i,j) = \Lambda_{1}(y,y') +  \\
& \quad2\kappa_{1}\kappa_{2}
\sum_{1\leq i<j\leq k}
(\sin\alpha_{i}\sin\alpha_{j})^{2} t(i,i+k)t(j,j+k)\cos\alpha_{i}\cos\alpha_{j}.
 \end{split} 
\end{equation}
From (\ref{105}) and (\ref{106}) we deduce that
 \begin{equation}\label{107}
\begin{split}
 & \Lambda(y,y')  = \Lambda_{1}(y,y')\quad  + \\
& \quad
 \kappa_{1}\kappa_{2}\bigg[  \sum_{1\leq i \leq k}t(i,k+i)^{2}\bigg(1 - \prod_{1\leq j\leq k, j\neq i} \sin ^{2}\alpha_{j} \bigg)
\sin^{4}\alpha_{i}  + \\
& \sum_{1\leq i<j\leq k}(\sin\alpha_{i}\sin\alpha_{j})^{2} t(i,i+k)t(j,j+k)\cos\alpha_{i}\cos\alpha_{j}\bigg] \quad + \\
& \quad \sum_{1\leq i<j\leq k}\langle \tilde{t}(i,j),\tilde{t}(i,j)B(i,j)\rangle \times \bigg[(1- \kappa_{1} - \kappa_{2}
+  \kappa_{1}\kappa_{2}(\sin\alpha_{i}\sin\alpha_{j})^{2})
\quad - \\
& \qquad (1  - \kappa_{1} -\kappa_{2} + \kappa_{1}\kappa_{2}\prod_{i=1}^{k} \sin ^{2}\alpha_{i} )\bigg] \\
& = \Lambda_{1}(y,y')+\Lambda_{2}(y,y')+\Lambda_{3}(y,y')
\end{split}
\end{equation}
and the proof of formula (\ref{96}) is completed. \\
The fact (see Lemma \ref{L5}) that the matrices $B(i,j)$ are positive-definite implies 
straightforwardly that $\Lambda_{3}(y,y')\geq 0$. It is perhaps a little less obvious
that $\Lambda_{2}(y,y')\geq 0$. This is a direct consequence of the lemma which follows. This result is certainly known but not having a reference at hand we give a proof of it in an appendix.
\end{proof}
 
\begin{lemma}\label{L7}
With the choice 
$k\geq 3$
and $ 0< \alpha_{i}< \pi/2, i=1,\dots,k$, 
the symmetric matrix 
$$
M(k)=\big( m^{i}_{j}\big)_{i,j=1,\dots,k}
$$
with 
\begin{equation}\label{108}
m^{i}_{i}=\bigg(1 - \prod_{1\leq j\leq k, j\neq i} \sin ^{2}\alpha_{j} \bigg)
\quad and \quad m^{i}_{j}= \cos\alpha_{i}\cos\alpha_{j}, \quad for \quad i\neq j,
\end{equation}
is positive-definite.
\end{lemma}
\subsection{The general case.}\label{SS2}
As announced with respect to the formulas (\ref{I-29}) and (\ref{I-30})
those expressions are  scalar products of two vectors $\tilde{y},\tilde{y}' \in \mathbb{R}^{p}$. The main point which allows to extend results of \ref{SS1} to the general case  is  the fact  that $\tilde{y}$ and $\tilde{y}'$ may be
expressed themselves via vectors $\tilde{y}_{0}$ and $\tilde{y}_{0}'$  given by formulas (\ref{62})-(\ref{64}). \\
\textbf{Consider at first the Case I (\ref{I-22}).}
\begin{lemma}\label{L8}
With the notations (\ref{I-25}) and (\ref{I-26}) the expression (\ref{I-29}) coincides with  the scalar product of the vectors
\begin{equation}\label{109}
\begin{split}
& \tilde{y}=\big(\tilde{y_{0}},c,d\sqrt{1-\kappa_{1}},e\sqrt{1-\kappa_{1}-
\kappa_{2}}\big)\\
& \qquad and\\
& \tilde{y}'=\big(\tilde{y_{0}}',c',d'\sqrt{1-\kappa_{1}},e'\sqrt{1-\kappa_{1}-
\kappa_{2}}\big)
\end{split}
\end{equation}
where $\tilde{y}_{0}$ and $\tilde{y}_{0}'$ are given by formulas (\ref{62})-(\ref{64})
of lemma \ref{L3}.
\end{lemma}
\begin{proof}
The proof lies on the obvious fact  
that the multivectors 
\begin{equation}\label{110}
\begin{split}
& u_{2m}\wedge\bigwedge_{i=1}^{k}u_{1i}\bigwedge_{i=1}^{k_{3}}u_{3i}
 \quad m=1,\dots,k \quad  and \\
 &  u_{lm}\wedge\bigwedge_{i=1}^{k}u_{1i}\bigwedge_{i=1}^{k_{3}}u_{3i}
 \quad l=4,5, \quad m=1,\dots,k_{l}
 \end{split}
\end{equation}
are mutually orthonormal.
This implies that
 \begin{equation}\label{111}
\begin{split}
& \langle y\wedge\bigwedge_{i=1}^{k}u_{1i}\bigwedge_{i=1}^{k_{3}}u_{3i},
y'\wedge\bigwedge_{i=1}^{k}u_{1i}\bigwedge_{i=1}^{k_{3}}u_{3i}\rangle \\ 
& =\langle y_{0}\wedge\bigwedge_{i=1}^{k}u_{1i}\bigwedge_{i=1}^{k_{3}}u_{3i},
y_{0}'\wedge\bigwedge_{i=1}^{k}u_{1i}\bigwedge_{i=1}^{k_{3}}u_{3i}\rangle
+ \langle d,d'\rangle + \langle e,e'\rangle\\
& = 
\langle y_{0}\wedge\bigwedge_{i=1}^{k}u_{1i},
y_{0}'\wedge\bigwedge_{i=1}^{k}u_{1i}\rangle
+ \langle d,d'\rangle + \langle e,e'\rangle
\end{split}
\end{equation}
and on the same basis
\begin{equation}\label{112}
\begin{split}
& \langle y\wedge\bigwedge_{i=1}^{k}v_{i}\bigwedge_{i=1}^{k_{3}}u_{3i}\bigwedge_{i=1}^{k_{4}}u_{4i},  y'\wedge\bigwedge_{i=1}^{k}v_{i}\bigwedge_{i=1}^{k_{3}}u_{3i}\bigwedge_{i=1}^{k_{4}}u_{4i}\rangle\\
& =\langle y_{0}\wedge\bigwedge_{i=1}^{k}v_{i},
y_{0}'\wedge\bigwedge_{i=1}^{k}v_{i}\rangle
 + \langle e,e'\rangle.
\end{split}
\end{equation}
From (\ref{I-25}), (\ref{I-26}) we obtain
\begin{equation}\label{113}
\langle y,y'\rangle = \langle y_{0},y_{0}'\rangle + \langle c,c'\rangle + \langle d,d'\rangle + \langle e,e'\rangle.
\end{equation}
Therefore 
\begin{equation}\label{114}
\begin{split}
 \langle y,y'\rangle & - 
\kappa_{1}\langle y\wedge\bigwedge_{i=1}^{k}u_{1i}\bigwedge_{i=1}^{k_{3}}u_{3i},
y'\wedge\bigwedge_{i=1}^{k}u_{1i}\bigwedge_{i=1}^{k_{3}}u_{3i}\rangle\\
& - \kappa_{2}\langle y\wedge\bigwedge_{i=1}^{k}v_{i}\bigwedge_{i=1}^{k_{3}}u_{3i}\bigwedge_{i=1}^{k_{4}}u_{4i},  y'\wedge\bigwedge_{i=1}^{k}v_{i}\bigwedge_{i=1}^{k_{3}}u_{3i}\bigwedge_{i=1}^{k_{4}}u_{4i}\rangle\\
& = \langle \tilde{y_{0}},\tilde{y_{0}}'\rangle + \langle c,c'\rangle + 
(1-\kappa_{1})\langle d,d'\rangle + (1-\kappa_{1}- \kappa_{2
})\langle e,e'\rangle
\end{split} 
\end{equation}
from which (\ref{109}) follows.
\end{proof}
As a consequence of (\ref{109}) and (\ref{I-31}) a simple calculation, which we do not detail, gives the exact value of  $ \Lambda_{g1}(y,y')$.
\begin{lemma}\label{L9}
We have
\begin{equation}\label{115}
\begin{split}
& \Lambda_{g1}(y,y')=\Lambda(y_{0},y_{0}') \\
& + (\kappa_{1}+\kappa_{2})\parallel c\wedge c'\parallel^{2} +
\kappa_{2}(1-\kappa_{1})\parallel d\wedge d'\parallel^{2}
\\
& +\kappa_{2}\sum_{1\leq i\leq k_{3}}\sum_{1\leq j\leq k_{4}}
(c_{i}d_{j}' - c_{i}'d_{j})^{2} \\
& +m(a,b,c)+ m(a,b,d)
\end{split}
\end{equation}
where
\begin{equation}\label{116}
\begin{split}
& m(a,b,c)=\sum_{1\leq i\leq k}\sum_{1\leq j\leq k_{3}}\Big[
(\kappa_{1}+\kappa_{2}\cos^{2}\alpha_{i})(a_{i}c_{j}' - a_{i}'c_{j})^{2}\\
 & \qquad + (\kappa_{2}+\kappa_{1}\cos^{2}\alpha_{i})(b_{i}c_{j}' - b_{i}'c_{j})^{2}\\
& \qquad + 2(\kappa_{1}+\kappa_{2})(a_{i}c_{j}' - a_{i}'c_{j})(b_{i}c_{j}' - b_{i}'c_{j})\cos\alpha_{i}\Big]\geq 0,\\
& m(a,b,d) = \kappa_{2}\sum_{1\leq i\leq k}\sum_{1\leq j\leq k_{4}}\Big[
(1-(1-\kappa_{1})\sin_{2}\alpha_{i})(a_{i}d_{j}' - a_{i}'cd_{j})^{2}\\
& \qquad + (1-\kappa_{1}\sin^{2}\alpha_{i})(b_{i}d_{j}' - b_{i}'cd_{j})^{2}\\
& \qquad + 2(a_{i}d_{j}' - a_{i}'d_{j})(b_{i}d_{j}' - b_{i}'d_{j})\cos\alpha_{i}\Big]\geq 0
\end{split}
\end{equation}
and $\Lambda(y_{0},y_{0}')\geq0$ is given by (\ref{95}) and (\ref{96}).
\end{lemma}
\textbf{The Case II (\ref{I-23}).}

\begin{lemma}\label{L10}
The expression (\ref{I-30}) coincides wwith  the scalar product of the vectors
\begin{equation}\label{117}
\begin{split}
& \tilde{y}=(\tilde{y_{0}},\sqrt{1-\kappa_{1}}d,e\sqrt{1-\kappa_{1}-
\kappa_{2}+  \kappa_{1}\kappa_{2}\prod_{i=1}^{k}(\sin\alpha_{i})^{2}})\\
& \qquad and\\
& \tilde{y}'=(\tilde{y_{0}}',\sqrt{1-\kappa_{1}}d',e'\sqrt{1-\kappa_{1}- \kappa_{2} + \kappa_{1}\kappa_{2}\prod_{i=1}^{k}(\sin\alpha_{i})^{2}}
\end{split}
\end{equation}
where $\tilde{y}_{0}$ and $\tilde{y}_{0}'$ are given by formulas (\ref{62})-(\ref{64})
of lemma \ref{L3}.
\end{lemma}
\begin{proof}
In this case we have by (\ref{I-27}) and (\ref{I-28})
\begin{equation}\label{118}
\langle y,y'\rangle = \langle y_{0},y_{0}'\rangle   + \langle d,d'\rangle + \langle e,e'\rangle
\end{equation}
and straightforwardly:
\begin{equation}\label{119}
\langle y\wedge\bigwedge_{i=1}^{k}u_{1i} ,
y'\wedge\bigwedge_{i=1}^{k}u_{1i}\rangle   
 = 
\langle y_{0}\wedge\bigwedge_{i=1}^{k}u_{1i},
y_{0}'\wedge\bigwedge_{i=1}^{k}u_{1i}\rangle
+ \langle d,d'\rangle + \langle e,e'\rangle,
 \end{equation}
\begin{equation}\label{120}
  \langle y\wedge\bigwedge_{i=1}^{k}v_{i} \bigwedge_{i=1}^{k_{4}}u_{4i},  y'\wedge\bigwedge_{i=1}^{k}v_{i}
\bigwedge_{i=1}^{k_{4}}u_{4i}\rangle
 =\langle y_{0}\wedge\bigwedge_{i=1}^{k}v_{i} ,  y_{0}'\wedge\bigwedge_{i=1}^{k}v_{i}\rangle  + \langle e,e'\rangle
\end{equation}
and
\begin{equation}\label{121}
\begin{split}
& \langle y\wedge\bigwedge_{i=1}^{k}u_{1i}\bigwedge_{i=1}^{k}v_{i}\bigwedge_{i=1}^{k_{4}}u_{4i},
y'\wedge\bigwedge_{i=1}^{k}u_{1i}\bigwedge_{i=1}^{k}v_{i}\bigwedge_{i=1}^{k_{4}}u_{4i}\rangle \\
& = 
\langle y_{0}\wedge\bigwedge_{i=1}^{k}u_{1i}\bigwedge_{i=1}^{k}v_{i} ,
y_{0}'\wedge\bigwedge_{i=1}^{k}u_{1i}\bigwedge_{i=1}^{k}v_{i} \rangle  
+ \langle e,e'\rangle.
\end{split}
\end{equation}
Accordingly, formulas (\ref{117}) follow directly from (\ref{I-30}) and (\ref{118}) - (\ref{121}).
\end{proof}
From (\ref{117}) and (\ref{I-32}) we get, by an easy calculus which we de not detail either.
\begin{lemma}\label{L11}
We have
\begin{equation}\label{122}
\begin{split}
& \Lambda_{g2}(y,y')=\Lambda(y_{0},y_{0}')  +  
\kappa_{2}(1-\kappa_{1})(1-\kappa_{1}\prod_{j=1}^{k}\sin^{2}\alpha_{j})\parallel d\wedge d'\parallel^{2}
\\
& \qquad + m(a,b,d)
\end{split}
\end{equation}

where
\begin{equation}\label{123}
\begin{split}
& m(a,b,d) = \kappa_{2}\sum_{1\leq i\leq k}\sum_{1\leq j\leq k_{4}}\Big[
\Big(\cos^{2}\alpha_{i} + \kappa_{1}(\sin^{2}\alpha_{i}- \prod_{j=1}^{k}\sin^{2}\alpha_{j})\Big)(a_{i}d_{j}' - a_{i}'cd_{j})^{2}\\
& \qquad + (1-\kappa_{1}\sin^{2}\alpha_{i})(1-\kappa_{1}\prod_{j=1}^{k}\sin^{2}\alpha_{j})(b_{i}d_{j}' - b_{i}'cd_{j})^{2}\\
& \qquad + 2(1-\kappa_{1}\prod_{j=1}^{k}\sin^{2}\alpha_{j})(a_{i}d_{j}' - a_{i}'d_{j})(b_{i}d_{j}' - b_{i}'d_{j})\cos\alpha_{i}\Big]\geq 0
\end{split}
\end{equation}
and $\Lambda(y_{0},y_{0}')\geq0$ is given by (\ref{95}) and (\ref{96}).
\end{lemma}
\section{Concluding remarks}
Following the argument developed in \cite{ref12} (section $5$) we can see that the inequality (\ref{I-8})
is still valid (for $n=2$ and $n=3$) in the setting of general discrete determinantal processes (finite or infinite)), namely   to discrete determinantal processes  associated to  positive contractions. It remains to investigate for the basic determinantal processes the case $n\ge4$ which does not seems very easy to handle.

\section*{Appendix: proof of lemma \ref{L7}}
The proof proceeds by induction. First, observe that if in a real symmetric positive-definite matrix we replace the diagonal entries respectively by greater elements then the new matrix obtained in this way remains positive-definite.
Now if in a M(k) matrix we delete the first row and the first column then 
the resulting matrix coicides with the
$M(k-1)=\big( m^{i}_{j}\big)_{i,j=2,\dots,k}$ matrix in which the diagonal entries 
$$
m^{i}_{i}=\bigg(1 - \prod_{2\leq j\leq k, j\neq i} \sin ^{2}\alpha_{j} \bigg), \quad
i=2,\dots,k
$$
are replaced respectively  by  the greater elements $\bigg(1 - \sin^{2}\alpha_{1}\prod_{2\leq j\leq k, j\neq i} \sin ^{2}\alpha_{j} \bigg)$.
Consequently, if the $M(k-1)$ matrix is positive-definite then the leading principal minors of $M(k)$ of orders less or equal to k-1 are positive and thus in order to prove that $M(k)$ is positive-definite it is enough to shown that the determinant $\det(M(k)$
is positive.
Now, for $k=3$ an elementary computation gives 
\begin{equation}\label{a2}
\begin{split}
\det M(3)=&
\bigg(\prod_{i=1}^{3} \sin ^{2}\alpha_{i}  + 
(\sin \alpha_{1}\sin \alpha_{2})^{2} + (\sin \alpha_{1}\sin \alpha_{3})^{2}
+ (\sin \alpha_{2}\sin \alpha_{3})^{2}
\bigg)
\\
& \times \prod_{i=1}^{3} \cos ^{2}\alpha_{i} >0 .
\end{split} 
\end{equation}
and it is obvious that the leading principal minors of $M(3)$ are positive.

To see that $\det M(k)>0$ for $k>3$, denote $a_{i}=\cos^{2}\alpha_{i}$ and
$b_{i}=\sin^{2}\alpha_{i}$, $i=1,\dots,k$, and suppose without loss of generality, permuting rows and columns (with the same index) of $M(k)$, if necessary,  that 
\begin{equation}\label{a3}
1>b_{1}\geq b_{2}\geq \dots \geq b_{k}>0.
\end{equation}
We obtain then
$$
\det M(k)=\prod_{i=1}^{k}a_{i}\times \det M_{0}(k)
$$
where
$$
 M_{0}(k)=\big( n^{i}_{j}\big)_{i,j=1,\dots,k}
$$
such that 
$$
n^{i}_{i}=\bigg(1 - \prod_{1\leq j\leq k, j\neq i} b_{j} \bigg)/a_{i}
  \quad and \quad n^{i}_{j}= 1, \quad for \quad i\neq j.
$$
It is well known and easy to see (for example by Sherman-Morrison formula) that
$$
\det M_{0}(k)=  \left[ n^{k}_{k}
+ (n^{k}_{k}-1)\sum_{i=1}^{k-1}\dfrac{1}{n^{i}_{i}-1}\right]
\times \prod_{i=1}^{k-1}(n^{i}_{i}-1).
$$
Moreover, the inequalities  (\ref{a3}) imply that
for $i=1,\dots,k-1$ we have
$$
n^{i}_{i}-1=\dfrac{b_{i}-\prod_{1\leq j\leq k, j\neq i} b_{j}}{1-b_{i}}>0.
$$
Hence, to finish the proof it suffices
 to shown that
$$
f(b_{k})=1 - \prod_{i=1}^{k-1} b_{i} + (b_{k}- \prod_{i=1}^{k-1} b_{i})\sum_{i=1}^{k-1}\dfrac{1-b_{i}}{b_{i}-\prod_{1\leq j\leq k, j\neq i} b_{j}}>0.
$$
Fix $b_{1},\dots, b_{k-1}$  and observe that for all $i=1,\dots,k-1$  the function 
$x \rightarrow  \dfrac{ x- \prod_{i=1}^{k-1} b_{i}}{b_{i}-x\prod_{1\leq j\leq k-1, j\neq i} b_{j}}$, $0\leq x <b_{k-1}$ is increasing and consequently
\begin{equation}
 f(b_{k})\geq f(0)=
1 - \prod_{i=1}^{k-1} b_{i} - \prod_{i=1}^{k-1} b_{i}\times\sum_{i=1}^{k-1}\dfrac{1-b_{i}}{b_{i}}.
\end{equation}
Puting $0<c_{i}=\dfrac{1-b_{i}}{b_{i}}$ we obtain
\begin{equation} 
 f(0)=\prod_{i=1}^{k-1}(1+ c_{i})^{-1}\times\bigg[\prod_{i=1}^{k-1}(1+ c_{i}) - 1 - \sum_{i=1}^{k-1}c_{i}\bigg]>0
\end{equation}
as desired.

\end{document}